\documentclass[12pt,letterpaper]{article}
\usepackage[a4paper, total={7in, 10in}]{geometry}

\usepackage{graphicx}
\usepackage{helvet}
\usepackage{authblk}
\usepackage{amsmath} 
\usepackage{amssymb} 
\usepackage{orcidlink} 
\usepackage[super,comma,sort&compress]  
   {natbib}
\usepackage[right]{lineno} \linenumbers
\usepackage{tabularx}

\makeatletter
\renewcommand{\maketitle}{\bgroup\setlength{\parindent}{0pt}
\begin{flushleft}
  \textbf{\@title}
  \@author
\end{flushleft}\egroup}
\makeatother

\usepackage{tikz}
\usepackage{mathtools}
\usepackage{hyperref}

\newtheorem{definition}{Definition}

\usetikzlibrary{shapes.geometric}
\tikzstyle{startstop} = [rectangle, rounded corners, minimum width=2cm, minimum height=1cm,text centered, draw=black, text=white, fill=black!60!blue!100, text width=4cm]
\tikzstyle{startstopsmall} = [rectangle, rounded corners, minimum width=2cm, minimum height=1cm,text centered, draw=black, text=white, fill=black!50!blue!100, text width=2cm]
\tikzstyle{io} = [trapezium, trapezium left angle=70, trapezium right angle=110, minimum width=2cm, minimum , text=white,height=1cm, text centered, draw=black, fill=black!60!blue!100]
\tikzstyle{process} = [rectangle, minimum width=2cm, minimum height=1cm, text centered, text=white, draw=black, fill=black!60!blue!100, text width=3cm]
\tikzstyle{decision} = [diamond, minimum width=2cm, minimum height=1cm, text centered, text=white, draw=black, fill=black!60!blue!100, text width=1.5cm]
\tikzstyle{arrow} = [thick,->,>=stealth]

\title{Robust Capacity Expansion Modelling for Renewable Energy Systems \\}
\date{}

\author[1,*,C,\orcidlink{0000-0002-3788-0473}]{Sebastian Kebrich}
\author[2,**,C,\orcidlink{0009-0007-7705-4508}]{Felix Engelhardt}
\author[1]{David Franzmann}
\author[2]{Christina Büsing}
\author[1]{Jochen Linßen}
\author[1,3,\orcidlink{0000-0001-9536-9465}]{Heidi Heinrichs}

\affil[1]{Forschungszentrum Jülich GmbH, Institute of Climate and Energy Research, Jülich Systems Analysis (ICE--2), Jülich, Germany}
\affil[2]{RWTH Aachen University, Faculty of Computer Science, Teaching and Research Area Combinatorial Optimization, Aachen, Germany}
\affil[3]{University of Siegen, Department of Mechanical Engineering, Chair for Energy Systems Analysis, Siegen, Germany}

\affil[C]{These authors contributed equally}

\affil[*]{Correspondence: s.kebrich@fz-juelich.de, lead contact}
\affil[**]{Correspondence: engelhardt@combi.rwth-aachen.de}

\begin{document}

\maketitle

\section*{SUMMARY}
Future greenhouse gas neutral energy systems will be dominated by renewable energy technologies providing variable supply subject to uncertain weather conditions. For this setting, we propose an algorithm for capacity expansion planning: We evaluate solutions optimised on a single years' data under different input weather years, and iteratively modify solutions whenever supply gaps are detected. These modifications lead to solutions with sufficient capacities to overcome periods of cold dark lulls and seasonal demand/supply fluctuations. 
A computational study on a German energy system model for $40$ operating years shows that preventing supply gaps, i.e. finding a robust system, increases the total annual cost by $1.6-2.9\%$. In comparison, non-robust systems display loss of load close to $50\%$ of total demand during some periods.
Results underline the importance of assessing the feasibility of energy system models using atypical time-series, combining dark lull and cold period effects.
\section*{KEYWORDS}
Resilience, Robustness, Bilevel, Operations Research, Critical time periods, Optimization under Uncertainty, Security of Supply
\section*{INTRODUCTION}

Modelling energy systems requires access to different types of data. For future greenhouse gas neutral energy systems, time--series data is of special importance, as a vast expansion of wind power and \emph{solar photovoltaic (PV)} that are dependent on the weather is perceived as indispensable~\citep{Gielen2019,Victoria_2022}. 
However, we can not expect to have accurate long--term information on future weather conditions. Instead, historical time--series serve as a substitute in capacity expansion planning, and often a single ''representative'' sample year is selected for optimisation~\citep{RINGKJOB_2018}. For example, the \emph{International Renewable Energy Agency (IRENA)} recommends using $2018$ as a reference year because it represents generation from renewable technologies well on average~\citep{IRENA2018Refyear}. 
In the following, we argue that using ''representative'' sample years may lead to \emph{energy system models (ESMs)} that appear sound at first glance, but would fail to meet supply in reality.
We then propose a tractable approach to counteract this effect and to make energy systems robust against uncertain conditions during operation, i.e. ensure security of supply.

There is ample evidence in literature that the results of energy system optimisation are sensitive to changes in weather time--series data.
Schyska et al. (2021) evaluate the sensitivity of capacity expansion models with regards to multiple sources of uncertainty. They conclude that some years are unsuited as reference years, as using them for optimisation leads to significant misallocation of assets.\cite{schyska2021sensitivity} 
Ruggles et al. (2024) assess how many years of weather data are needed to ensure ESMs are reliable even out of sample, i.e. if a different weather year were to realise. They conclude that between $15$ and $40$ years are required depending on amount of imports available/the acceptable level of loss of load ~\cite{Ruggles2024ReliableWindSolar}.
The works by Haddeland et al. (2022)~\cite{haddeland2022effects}, and Staffel et al. (2018) ~\cite{staffell2018increasing}, which builds on earlier work by Pfenninger et al. (2017)~\cite{PFENNINGER20171},  similarly find the choice of weather years significantly effects renewable generation and thus power output. 
Not only that, but the effect of the weather increases with increasing share of renewable technologies~\cite{collins2018impacts}.
De Marco et al. (2025) identify energy shortages across Europe and use those to optimize climate--resilient energy systems stochastically~\cite{DeMarco2025ClimateResilient}.
Not all research agrees on this; Schlachtenberger et al.(2018)~\cite{schlachtberger2018cost} optimise three weather years with hourly data both individually and as one time--series with a resolution of 3h per time step, finding only small variations in \emph{total annual cost (TAC)} and installed capacities. However, they note that aggregating multiple hours together introduces a smoothing effect that systematically favours PV and underestimate battery and wind generation requirements. 

Additionally, there is a fundamental tension between identifying typical years, and build ESMs that are protected against extreme events, i.e. robust. On this, Hilbers et al. (2020)~\cite{Hilbers_2020} introduce a method of importance subsampling for time--series aggregation to explicitly preserve extreme events in the weather data as an alternative to established "representative days" clustering approaches. 
Ryberg's dissertation~\cite{Ryberg2020Diss} and Ryberg et al. (2019)~\citep{ryberg2019occurrence}investigate the impact of generation lulls in a energy system for a large part of Europe calculating backup capacities required to overcome these. Ruhnau et al. (2022)~\cite{Ruhnau_2022} look into the storage requirements for a renewable--based ESM for Germany using $35$ years of weather data taking consecutive extreme events into account. They conclude that consecutive extreme events increase storage requirements significantly compared to even the most extreme, but singular events. 
Thus, atypical weather time--series may be particularity well-suited for optimisation because they capture important system behaviours, e.g. dark lulls, that significantly impact ESMs.
In this context, Grochowicz et al. (2024)~\cite{grochowicz2023intersecting} discuss optimising sequential weather years. They use a geometry--based solution approach targeting the solution space. In  follow--up work, they use electricity shadow prices to identify difficult weather periods~\cite{grochowicz2024using}. They observe that such difficult weather periods are not just meteorological events, but results of the interplay of meteorology and electricity storage and network structures.

In stochastic optimisation, the notion of using typical, i.e expected, behaviour as a baseline for optimisation is a well known concept. 
Its usefulness is determined by the value of a stochastic solution, which represents the gap in expected performance of a solution obtained with expected data and one determined by solving a full stochastic optimisation problem ~\citep{Birge1982}.
In general, this gap may be arbitrarily large. However, research has also shown that we can still extract insights from solutions obtained with with such data ~\citep{Maggioni2013}.
Specifically concerning weather robustness, Forghani et al. (2025) propose an intermediate approach where not only a representative but also worst-/best-cost years are used as input data for a stochastic optimisation. Their results show that this reduces loss-of-load to practically zero, at $<1\%$ additional cost~\cite{forghani2025modelingrobustenergysystems}.

An important modelling decision in this context is the value and availability of recourse, i.e. the ability to specify or change parts of the solution if uncertainty realises. 
In capacity expansion planning, imports constitute such recourse. If arbitrary imports are allowed, any misplanning can be compensated for, if at a cost.
We focus on a different setting, the one where imports are limited or energy systems are to be self--sufficient.
Here, we use \emph{adaptive robust optimisation (\ref{def:aro})} to model limited recourse while still enforcing strict guarantees on security of supply.general solution approaches.
To the best of our knowledge, Zeyringer et al. (2018)~\cite{zeyringer2018designing} were the first to assess the effects of weather uncertainty in input data on capacity--expansion planning, looking at an ESM of Grein.They propose using multiple historical time series, optimising over each, and then evaluating the operational costs/supply costs incurred by the proposed energy systems under the assumption that the time--series of a different historical year realises. Then, an ESM with lowest worst--case costs/lowest supply gap amongst the solutions is selected. They find that starting with the ‘’wrong’’ reference year may lead to misallocation of resources and that modelling with multiple years leads to more consistent results with lower worst--case cost, the latter of which had also been noted before~\citep{Bloomfield_2016}.
For the US, Dowling et al. (2020)~\cite{Dowling_2020} make the case for multi--year modelling to accurately capture long--term storage effects, noting that the cost of variable renewable power systems are especially sensitive to long--duration storage costs. 
A recent work by Gøtske et al.~\cite{gøtske2024designingsectorcoupledeuropeanenergy} also assesses energy systems based on different weather years. They employ CO$_2$ emitting backup technologies, and analyse structural elements of the respective solutions.
Either approach allows to select more suitable ESMs, but it can not assure solutions meet certain supply/demand across all years.

\subsection*{The Two-Stage Robust Setting}
To illustrate the setting considered in this work, we begin with an example for an \ref{def:aro} ESM. \ref{def:aro}, alternatively called two-stage robustness, is characterised though a bilevel structure~\citep{Goerigk2025} in which a decision maker has to make a set of first stage decisions $x \in \mathcal{X}$, then one of multiple scenarios $u \in \mathcal{U}$ may realise and afterwards the decision maker has the option to react to scenario $u\in \mathcal{U}$ with a scenario-dependent set of decisions $y(u)\in \mathcal{Y}(x,u)$. Equation~\eqref{def:aro} illustrates the structure of a \ref{def:aro} problem. Here, $c$ and $d$ are vectors of first and second stage costs, respectively. The matrix $A$ encodes all information that is not affected by uncertainty, whereas the matrices $B(u)$ and the right-hand side vector $b(u)$ are dependent on the scenario $u$.
\begin{align}
    \min_{x \in \mathcal{X}} c^\top x \; +& \; \max_{u \in \mathcal{U}} \; \min_{y(u) \in \mathcal{Y}(x,u)} \; d^\top y(u) \label{def:aro}\tag{ARO}\\
    \text{s.t.} \quad & A x + B(u) y(u) \geq b(u)\; \forall u \in \mathcal{U}
    \notag
\end{align}
For the ESM(s) considered in this work, the first stage decisions correspond to the capacity expansion planning, e.g. investment in solar power plants, power lines, storage units, and the second stage corresponds to operational decisions, e.g. how much energy is to be stored in which battery at what time. 
If energy imports are allowed, these will also be encoded in $\mathcal{Y}$.
The scenario $u$ corresponds to a specific, continuous time period.
Table \ref{tab:aro_example} illustrates this setting. Here, we consider different locations $A, B, C$ for a solar power plant, whose production is subject to changing weather patterns during the time periods $u_1,u_2$. 

\begin{table}[h!]
    \centering
    \begin{tabular}{c|c|c|c}
       Scenario  & Location A & Location B & Location C \\ \hline
       $u_1$ & 200 & 75 & 50 \\
       $u_2$ & 50 & 75 & 150
    \end{tabular}
    \caption{Investment problem: Chose number of solar power plants to build on locations $A,B,C$. Values indicate expected energy supply in $\frac{\mathrm{GWh}}{\mathrm{time~period}}$ per investment.}
    \label{tab:aro_example}
\end{table}

The optimal solution depends on what type of recourse is available and how much supply is needed.
For example, if we assume similar cost at all locations and no/expensive imports, investing in location $B$ is the most cost effective for meeting a demand of $75\frac{\mathrm{GWh}}{\mathrm{time~period}}$. However, if the total demand were $200\frac{\mathrm{GWh}}{\mathrm{time~period}}$, it is more effective to invest in $A$ and $C$ instead of $B$. 

Returning to our example in Table \ref{tab:aro_example}, the approach by Zeyringer et al. would not be able to find the optimal robust solution for a total demand of $75$MWh/year, which is to invest in location $B$. This is because if the future is known, i.e. we assume either $u_1$ or $u_2$ to realise, the optimal investment is always to either invest in location $A$ or location $B$. 
An alternative approach was put forward by Caglayan et al. (2019)~\cite{Caglayan_2019impact}, who noted variations of up to $20\%$ variations in TAC, and significant differences in energy system designs when optimising an integrated electrical/hydrogen ESM for Western Europe with different weather years. They propose an iterative approach, where they begin with optimising individual years, but then average their designs and take the average as a lower bound for installed capacities in a new round of optimisation~\citep{Caglayan_2021robust}. Since the (average) capacities installed during sequential algorithm iterations will be non-decreasing and bounded, convergence and thus a feasible solution will be guaranteed in most cases. In our example, the algorithm would enforce capacity $>0$ for location A and C, which would eventually lead to a feasible but unnecessarily expensive solution that uses both A and C.

\subsection*{A Novel Algorithm for Computing Robust ESMs}
The proposed algorithm consists of four main steps, as shown in Figure \ref{fig:flowchart}. As its input, the algorithm takes a set of $n$ times series, one of which is designated as a starting/reference year. 
We begin with solving the capacity expansion planning problem for this year (CAPEX), then we evaluate the performance of the proposed investment in all $n$ years by solving the corresponding unit commitment problem (UC), which formalises the validation proposed in previous works ~\cite{zeyringer2018designing,gøtske2024designingsectorcoupledeuropeanenergy}. 

To ensure feasibility, existing approaches allow for energy imports. 
We consider the robust setting where imports are bounded/not available, i.e. the ESM needs to be self-sufficient. Based on the previous step, we identify supply gaps, which we denote by $\delta$. The underlying ESM for Germany for $2045$ utilizes renewable technologies only and includes PV, wind, Li--ion batteries, electrolyzer, hydrogen salt caverns, hydrogen combined cycle gas turbines as well as electricity transmission and hydrogen pipelines. The models are described in detail in the methods section for both \emph{gurobipy} and the ETHOS.FINE.
If the supply gaps are sufficiently small, we terminate. Otherwise, we use the information encoded in $\delta$ to iteratively make solutions more robust. 
For that, modifications (MOD) are applied to the optimisation problem by either adding artificial demand \ref{mod:1_demand}, substituting parts of the time series data \ref{mod:2_timeseries} or iteratively adding both demands and lazy valid inequalities~\ref{mod:3_combine}, before reoptimising CAPEX to evaluate the loss of load. Detailed definition of the modifications can be found in the methods section.

\begin{figure}[h!] 
\centering
\begin{tikzpicture}[node distance=1.5cm]
    \node (start) [startstop] {\small \textbf{Start}: Solve CAPEX for reference year};
    \node (proc2) [process, right of=start, xshift=3cm] {\small Solve UC for all $n$ years};
        \node (proc4) [decision, below of=proc2, yshift=-1.25cm] {\small Small loss of load?};
    \node (proc3) [process, right of=proc4, xshift=3cm] {\small Apply MOD1-3 to CAPEX$^*$};
    \node (proc5) [process, above of=proc3, yshift=1.25cm] {\small Reoptimise modified CAPEX$^*$};
    \node (stop) [startstop, left of=proc4, xshift=-3cm] {\small \textbf{End}: Robust energy system design};
    \draw [arrow] (proc4) -- node[anchor=north] {yes} (stop);   
    \draw [arrow] (proc4) -- node[anchor=north] {no} (proc3);   
    \draw [arrow] (proc3) -- (proc5);
    \draw [arrow] (proc2) -- (proc4);
    \draw [arrow] (proc5) -- (proc2);
    \draw [arrow] (start) -- (proc2);
\end{tikzpicture}
\caption{Flowchart depicting the proposed methodology for determining robust energy systems. In each main loop iteration, $n-1$ unit commitment problems (UC) and one modified capacity expansion problem (CAPEX$^*$) are solved until loss of load is sufficiently small.}
\label{fig:flowchart}
\end{figure}

The main contributions of our work include not only insights derived for the specific ESM, and an algorithm for solving CAPEX problems to find solutions that are robust to variations in operational conditions across the $n$ years, but foremost a systematic comparison of different approaches for designing robust ESMs. Based on this, we recommend three working modifications that all lead to robust energy systems systems.

\section*{RESULTS AND DISCUSSION}\label{sec:results}

In the following, we first compare solutions derived from individual years' data and assess their (lack of) robustness. Then, we modify them to become more robust, and evaluate the features of robust solutions.
Note that we provide results for two different models, since some modifications require custom callbacks/lazy constraints which are hard to integrate into modelling frameworks. Thus, we use one model directly coded in  \emph{gurobipy}  that allows for all modifications but is simplified, see Appendix \ref{app:lp_esm}, and a more realistic model implemented in the ETHOS.FINE framework.

Note that we also only use a single demand time-series from $2012$ for comparison. 
This significantly simplifies evaluation of the algorithms as otherwise, loss of load and installed would need to be considered relative to each year's demands.
However, in most practical settings where demands and weather are linked, e.g. when employing significant electric or hydrogen-based heating, the paired demand time-series ought to be used.

\subsection*{Optimal Capacities Strongly Depend on Yearly Weather Data}

The TAC for energy systems for the $38$ node model of Germany within ETHOS.FINE for the $40$ different years of time--series data deviates around an average of $106$bn€~with cost between $96.4$bn€~ and $113.6$bn€ annually, which equals $-9\%$ to $+7\%$ compared to the default reference year $2018$ recommended by IRENA~\cite{IRENA2018Refyear}.
The results of optimising each year independently are given by Figure \ref{fig:tac_40y}. While the variations in overall TAC are limited, the energy system designs show substantial deviations.
The cost shares of single technologies across the $40$ different single years vary by $69\%$ for hydrogen pipelines, $57\%$ for hydrogen salt caverns, $53\%$ for CCGT, $44\%$ for Li--ion batteries, $40\%$ for PV, $38\%$ for electrolysers, $22\%$ for the electricity grid, and $20\%$ for onshore wind, making it challenging to draw recommendations for planned capacity expansion for future energy systems. 

\begin{figure}[h!]
\includegraphics[width=\textwidth]{Figure_2.png}
\centering
  \caption{TAC comparison by technology for ESMs optimised from 1980--2019 aggregated for a $38$ node Germany model set up in ETHOS.FINE.}
  \label{fig:tac_40y}
\end{figure}

\begin{figure}[h!]
    \centering
    \includegraphics[width=\textwidth]{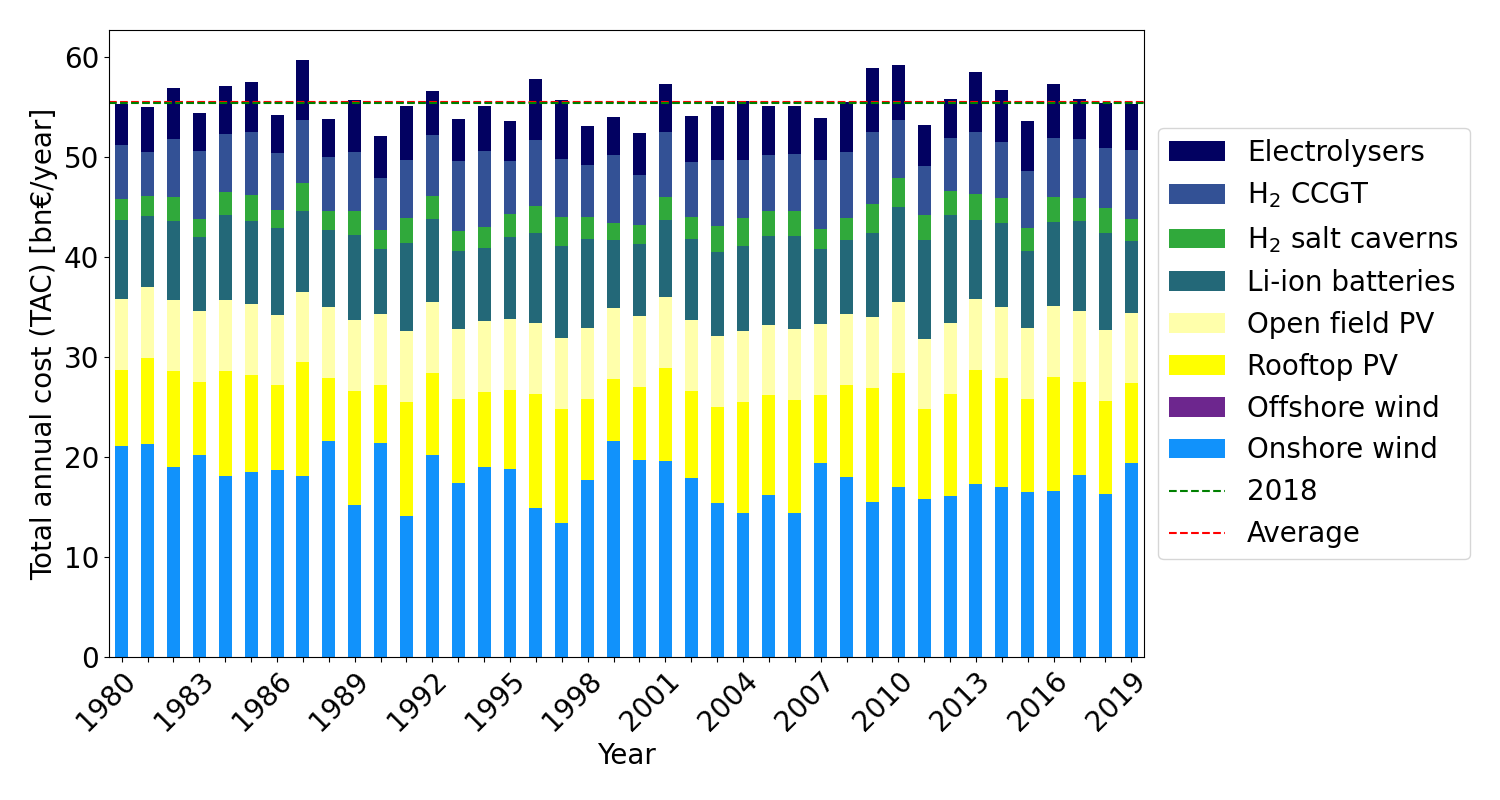}
    \caption{TAC comparison by technology for ESMs optimised from $1980-2019$ on a single node Germany model set up using gurobipy.}
    \label{fig:individual_years}
\end{figure}

In the simplified single node model, the cost variation for single technologies is on average slightly higher than in the 38-node model: onshore wind ($45\%$), rooftop PV ($57\%$), Li--ion batteries ($51\%$), hydrogen salt caverns ($54\%$), electrolysers ($54\%$) and CCGT ($53\%$). Open field PV is installed up to its maximum capacity for all years, while offshore wind is never utilised. 
A plot of all individual years is given in Figure~\ref{fig:individual_years}.

In summary, even small input data differences may lead to large variations in investment into different technologies. 
For example, offshore wind is not utilised at all in some years, but makes up over $13\%$ of the TAC in $2014$.
Note that the single node model is not directly comparable to the ETHOS.FINE model, as it contains various simplifications. As such, solutions are significantly cheaper averaging $55.4$bn€~with between $52.1$ and $59.7$bn€~annually, i.e. $-6\%$ to $+7\%$. 

\subsection*{Systems Optimised on Single Years Cause Supply Gaps}

\begin{figure}[h!]
    \centering
    \includegraphics[width=\textwidth]{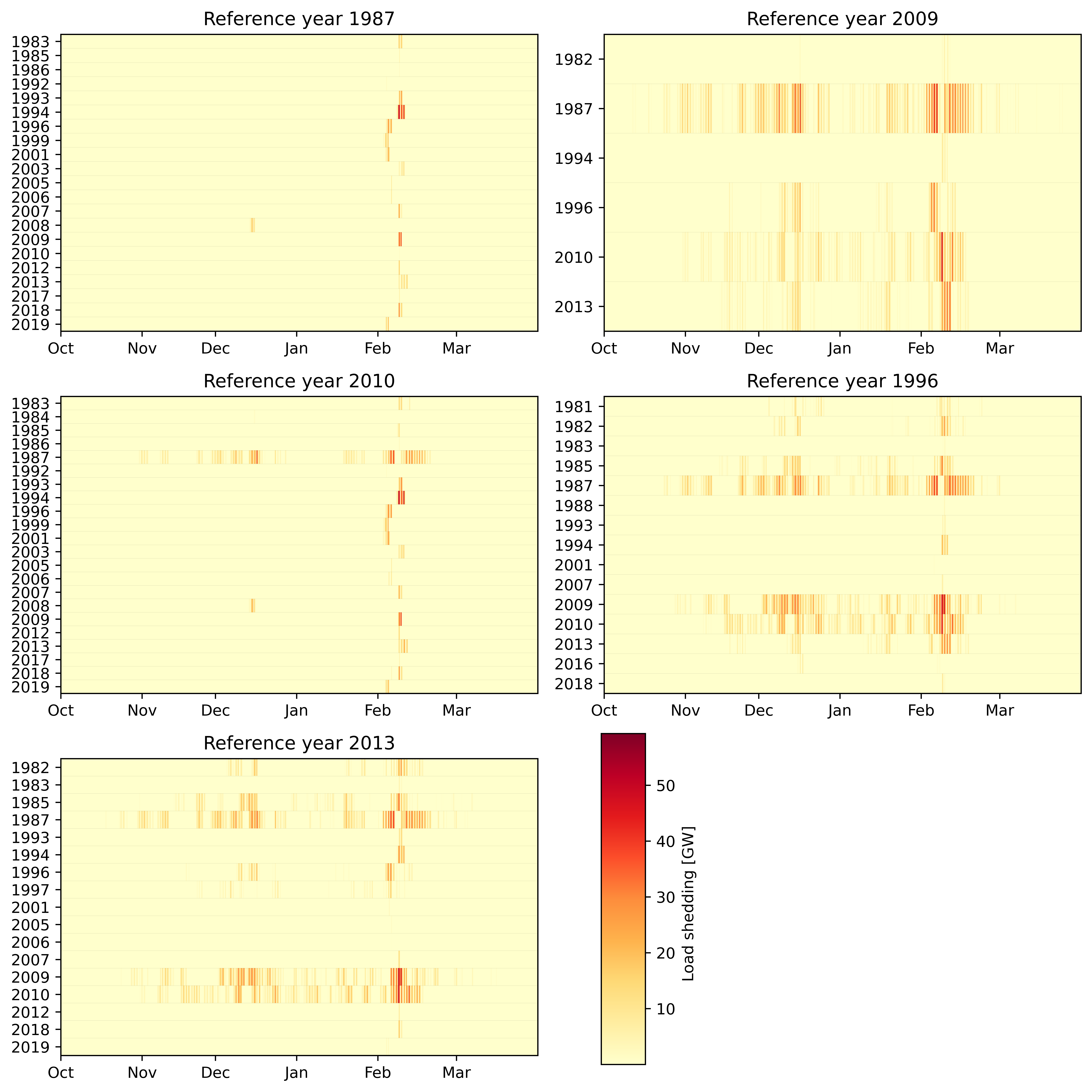}
    \caption{Load shedding when testing the feasibility of the $5$ selected reference years in ETHOS.FINE. Only the years and the months from October to March are included, where load shedding occurs. For comparison, hourly load in the demand data we use is roughly $100$GW.}
    \label{fig:Load_shedding}
\end{figure}

Feasibility testing shows that all $40$ energy system designs solely based on one year's time--series lead to supply gaps in multiple other years. This means that none of the system designs are robust. 
However, we are not only interested in whether there is any loss of load, but also in the magnitude of said loss of load. 
Figure \ref{fig:Load_shedding} shows the amount of load shedding when testing the feasibility of one of five reference years used in ETHOS.FINE during the other $39$ years.

Two trends can be observed: On the one hand, when the feasibility is tested for $1987$ and $2010$, load shedding occurs mostly in the first half of February. This is caused by demand spikes in early February due to the low temperatures in that time of the year, which the electricity demand was based on.
At the same time, the energy systems optimised for $1987$ and $2010$ are characterised by the highest amounts of PV and wind capacities installed while having insufficient backup capacities.

On the other hand, the energy system of the reference years $2009$ has increased backup capacities, therefore suffering less during cold darklull periods, but does not have enough PV and wind capacities to produce enough hydrogen and therefore has load shedding more evenly distributed over the months from October to March. The other two reference years $1996$ and $2013$ are in the middle of these two case.

\begin{figure}[h!]
\centering
  \includegraphics[width=\textwidth]{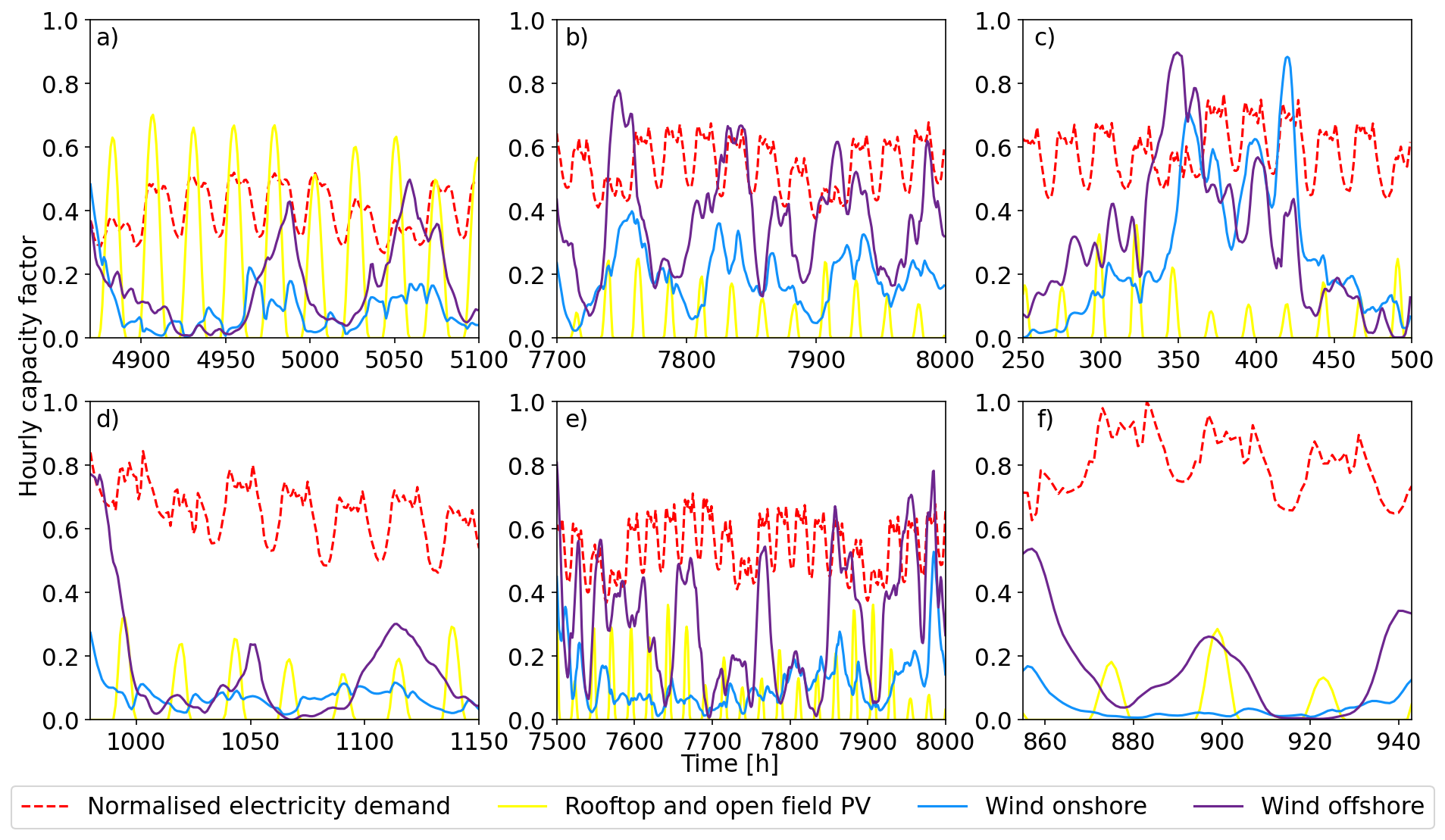}
  \caption{Six time periods from the $40$ weather years and $1$ year of future electricity demand data for Germany. The electricity demand is normalised to prevent overweighing; the weather data is aggregated. The upper three diagrams represent non--critical time periods, the lower three critical ones. 
  Subfigure a) is a typical summer period with high PV availability and low electricity demand due to low heating requirements. 
  In b) and c), typical autumn and winter period are shown. They are characterised by low availability of PV, but ample wind power to supply sufficient electricity. Note the increased electricity demand due to increased heating required. 
  In d), a typical dark lull  is characterised by low availability of PV and negligible amounts of wind, which coincides with high electricity demand due to increased heating. 
  Subfigure e) shows an elongated dark lull period. Low availability of both PV and wind combined with increased electricity demand lead to overall difficult period requiring large amounts of hydrogen to be burned in the energy system. The last graphic f) shows the most critical period in the 40 years of data. Negligible wind combined with low availability of PV coincide with the highest electricity demand in the data due to high heating demand during an extreme cold spell in all of Germany.}
  \label{fig:ExampleTimePeriods}
\end{figure}

Figure \ref{fig:ExampleTimePeriods} illustrates the reasons for that by showing three uncritical (top row) periods with sufficient supply as well as three critical (bottom row) time periods with significant supply gaps. 
During critical time periods electricity generation from PV and Wind is low, providing less than 50\% of electricity demand. The remaining electricity demand needs to be fulfilled by backup capacities, i.e. $H_2$ CCGT.
The periods have varying duration and they are identified via clustering and feasibility testing. 
Uncritical time periods are characterised by high availability of PV or onshore wind or both, while critical ones are characterised by low availability of PV and low to negligible onshore wind combined with high demand. Offshore wind plays only a minor role due to its limited utilisation. The most critical time period takes place in $1994$, see Figure \ref{fig:ExampleTimePeriods}f). Here, wind and solar supply indicate a dark lull. Combined with the cold period identified in the electricity data, this constitutes a cold dark lull.

For the  \emph{gurobipy}  implementation, we get similar results: The minimal annual supply gap across all years is $1.170$GWh for $1987$, which is about $0.13\%$ of the total annual electricity demand. Additionally, during shorter time-periods the supply gap can reach up to half of the required electricity demand.
In comparison, optimising using $1990$'s weather data leads to a peak supply gap of $69.500$GWh, i.e. $8\%$, if the weather of $1987$ were to realise.

In summary, this strongly supports the results by Cagalayan et al.~\cite{Caglayan_2019impact,Caglayan_2021robust} -- weather patterns matter and ignoring in energy systems with significant renewables may lead to systems that have large supply gaps under anything but the most optimal conditions, which we cannot simply assume to be covered through imports.

\subsection*{Effect of Modifications on Installed Capacities} 
All modifications lead to increased investment in CCGT, salt caverns and electrolysers, although the latter is less pronounced for \ref{mod:1_demand}, as seen in Table \ref{tab:all_algorithms}.
This is expected, as CCGT can provide electricity when the generation is low and demand is high and installing more CCGT also requires more H$_2$ infrastructure such as salt caverns for storage and electrolysers for H$_2$ conversion. 

Note that the more detailed  model implemented in ETHOS.FINE was used for evaluation of \ref{mod:2_timeseries}, since it outperformed the other algorithms when implemented in \emph{gurobipy} and is therefore viewed as the most promising approach. An additional advantage of \ref{mod:2_timeseries}: It can be easily integrated into existing modelling frameworks.

\paragraph{ \ref{mod:1_demand}} leads to robust energy systems regardless of the initial time--series chosen. 
On average, this incurs additional cost of $7.8$bn€ if no smoothing is performed.
The average total cost of a robust system reaches $63.1$bn€, with a range of $[60.9,69.7]$bn€.
In comparison, using smoothing leads to slightly more expensive solutions. On average, making an ESM robust incurs additional cost of $9.7$bn€ (+$17\%$). The average total cost of a robust system reaches $65.2$bn€, with a range of $[61.0, 72.7]$bn€.

\begin{figure}[h!]
\centering
  \includegraphics[width=\textwidth]{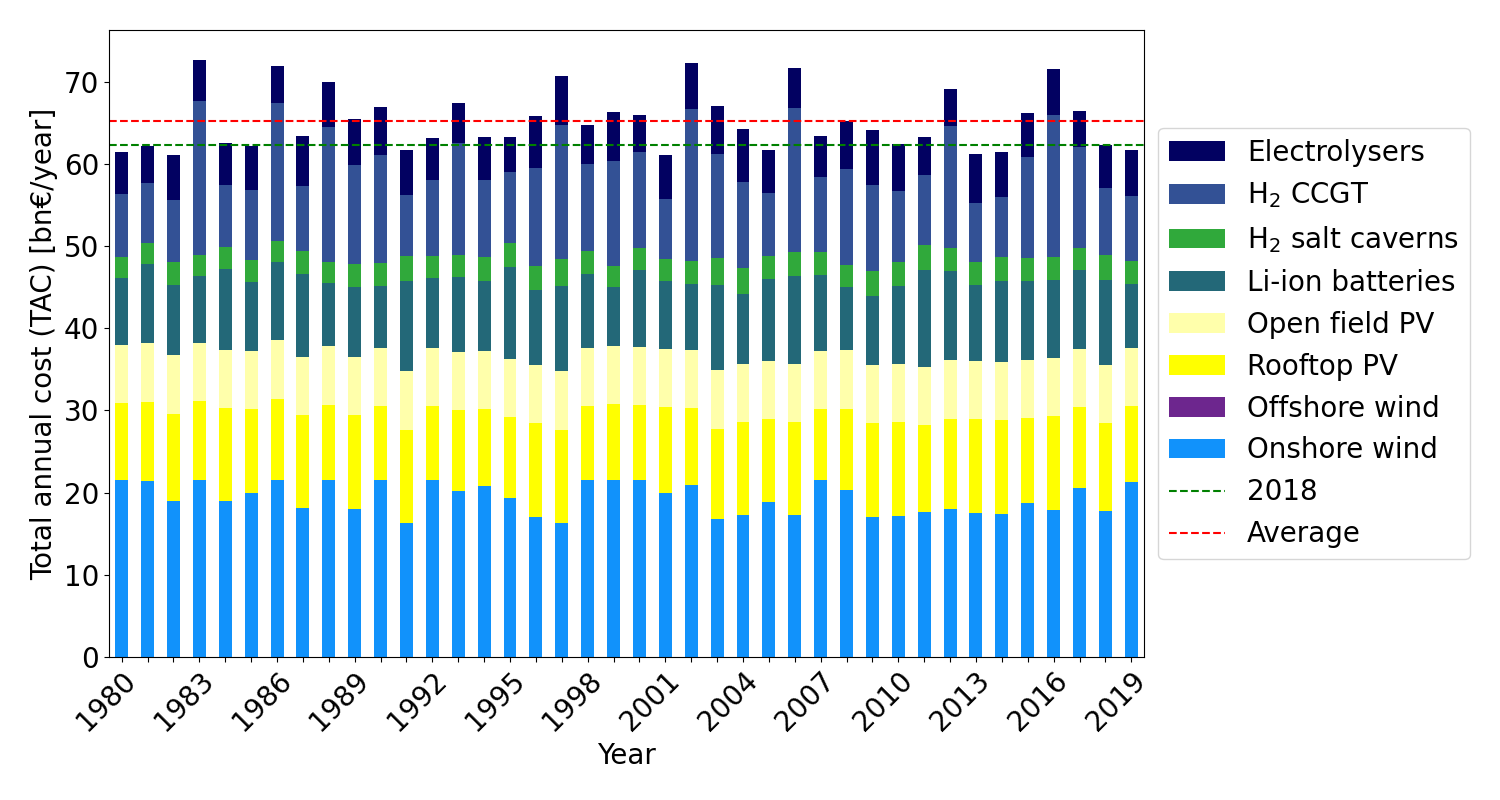}
  \caption{TAC comparison from $1980-2019$ for robust solutions using modification smoothed \ref{mod:1_demand} for the single node model in gurobipy, no temporal aggregation.}
  \label{fig:Alg2}
\end{figure}

Figure \ref{fig:Alg2} gives the results of optimising each year independently for the smoothed \ref{mod:1_demand}. 
While on average more expansive than non--smoothed  \ref{mod:1_demand}, the smoothed algorithm leads to on average lower investment increase in short--term battery storage ($+25\%$ non--smoothed  vs. $+12\%$ for smoothed) and a significantly higher investment in CCGT ($+86\%$ for non--smoothed vs. $+36\%$ for smoothed).
The higher costs may be due to the fact that additional artificial demand is added in time periods adjacent to those with previous supply gaps, which generates small supply demand gap time periods. 
The strong invest in CCGT compared to the non--smoothed modification suggests that CCGT power plants are used to offset those artificial demand gap time periods.

In summary, smoothed \ref{mod:1_demand} leads to feasible solutions, but care should be taken that artificial energy demands do not lead to excessive building of CCGT power plants.
Convergence of \ref{mod:1_demand} can be slow, sometimes taking more than $20$ iterations for a single pair of years. This number of iterations was used as a cut-off criterion, as no improvements were observed after that during testing. This appears to be due to very small residual supply gaps of a few $GWh$ that get found and added to the model repeatedly. Given the small size of those supply gaps, the large overall production, and the fact that Gurobi~ \cite{gurobi} was used as a solver, which does not perform exact arithmetic, this may be caused by numerical instabilities. Using a suitable termination criterion (e.g. number of iterations or total supply gap less than some small number of GWh) counteracts this.

Notably, non-smoothed  \ref{mod:1_demand} incurs a bias towards installing more Li--ion battery storage. This is to expected, as artificial short term demand peaks are added, and Li--ion batteries are well--suited to compensate for those.
Their capacity was increased by on average more than $25\%$, with a range of $[7.7,19.1]$bn€, compared to $[6.5,10.8]$bn€~in the reference years.
Finally, non-smoothed  \ref{mod:1_demand} finds the overall cheapest robust solution. That solution is characterised by slightly more investment in onshore wind capacity ($19.8$bn€, $+12\%$) and roof top PV ($10.2$bn€, $+10\%$) than in an average single year solution. No additional batteries are installed, but more electrolysers ($5.6$bn€, $+18\%$), CGGT ($7.3$bn€, $+23\%$) and salt caverns ($2.7$bn€, $+16\%$).

\begin{figure}[h!]
\includegraphics[width=\textwidth]{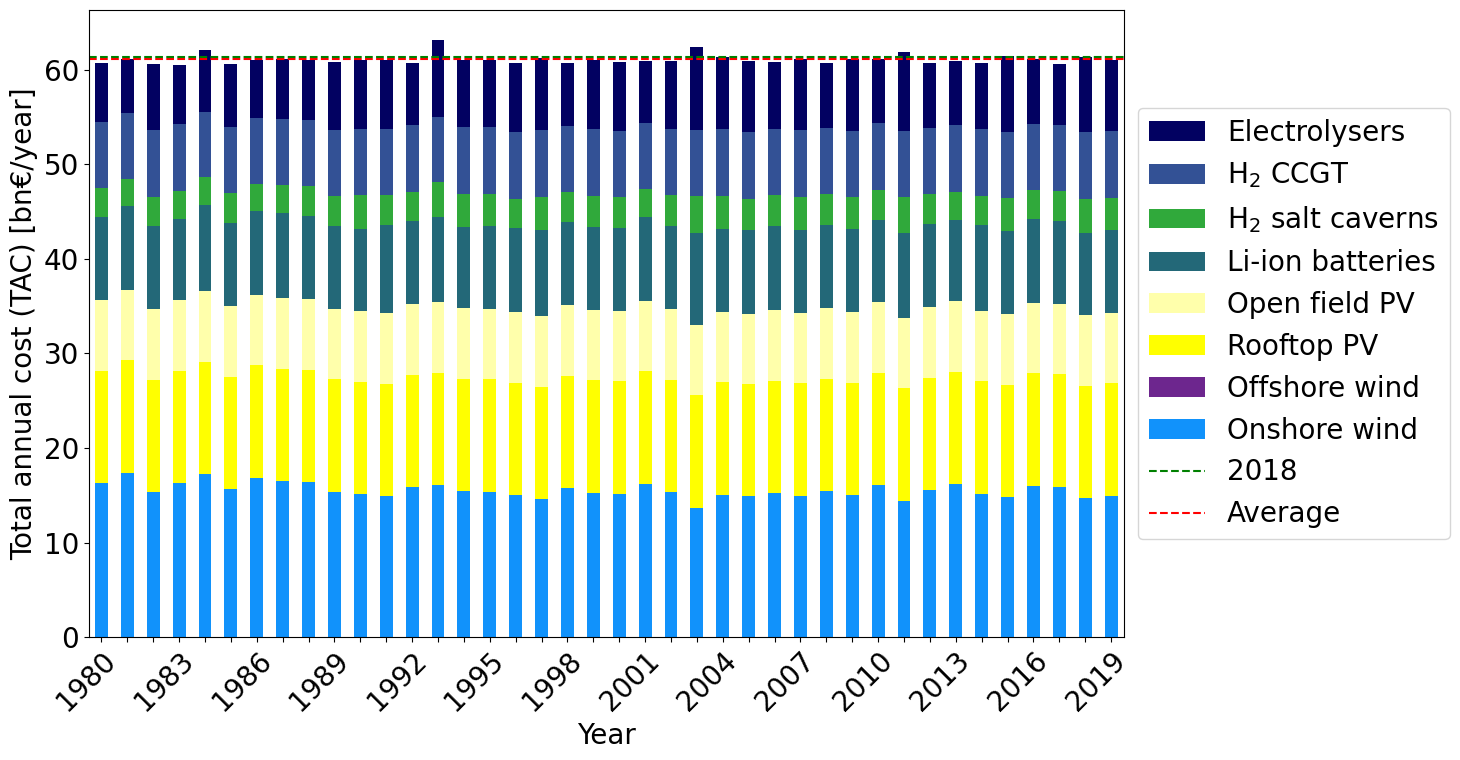}
\caption{TAC comparison for optimising for from $1980-2019$ for robust solutions using modification \ref{mod:2_timeseries} \emph{gurobipy} for the singe node model in gurobipy, no temporal aggregation.}
\label{fig:mod2_Results_gurobi}
\end{figure}

\begin{figure}[h!]
\includegraphics[width=\textwidth]{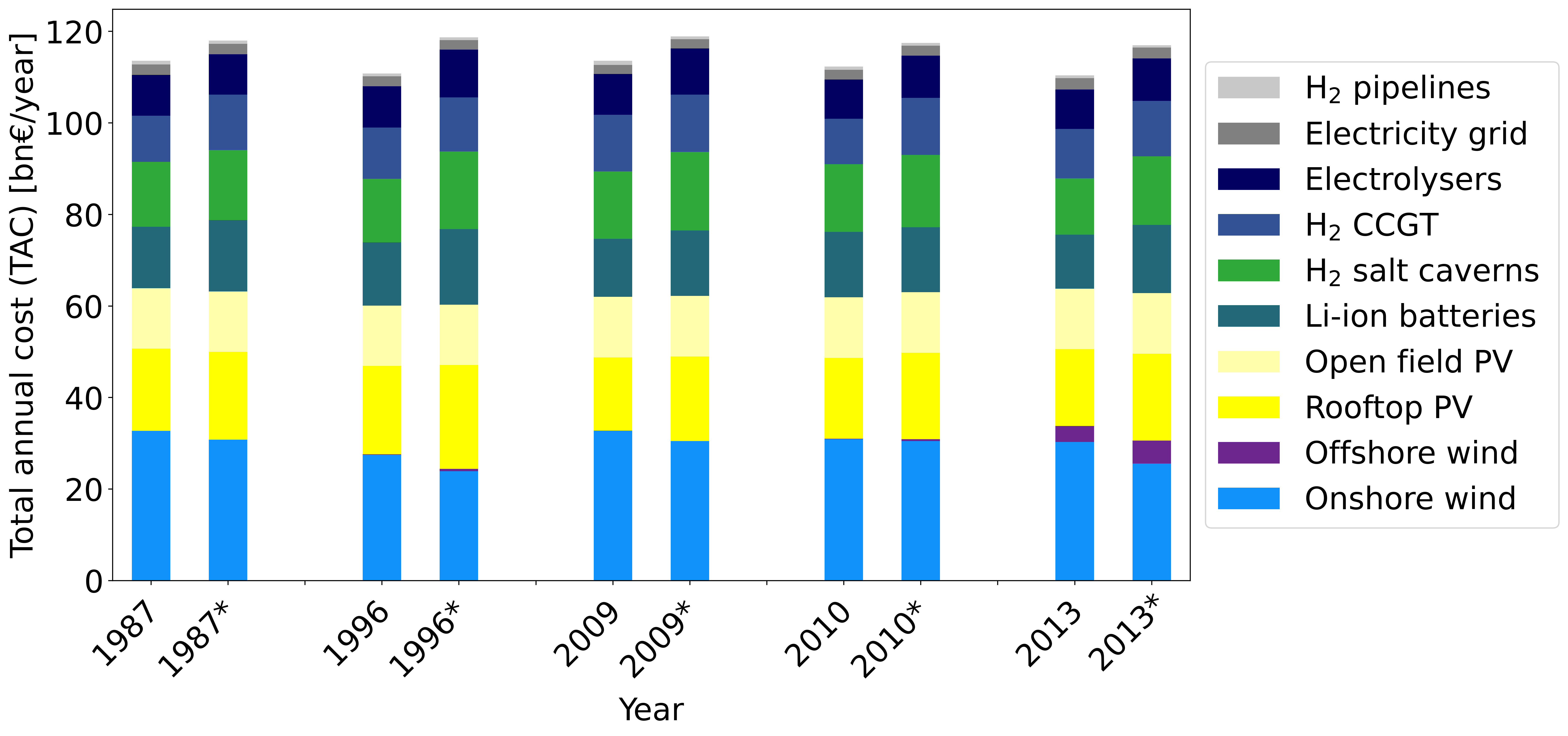}
\caption{TAC comparison for optimising for $5$ selected reference years for Germany. Columns marked with a * indicate the robust system designed via \ref{mod:2_timeseries} ETHOS.FINE.}
\label{fig:mod2_Results_FINE}
\end{figure}

\paragraph{\ref{mod:2_timeseries}} is evaluated in Figures \ref{fig:mod2_Results_gurobi} and \ref{fig:mod2_Results_FINE} for the model in  \emph{gurobipy}  and the model in ETHOS.FINE, respectively, and Table \ref{tab:Alg3_Results_capacities}, which give an overview of the results of generating robust ESMs.
In the ETHOS.FINE model the five most expensive solutions based on the single years $1987, 1996, 2009, 2010$ and $2013$ were selected of the $40$ energy system designs to make them robust using \ref{mod:2_timeseries}. These are hereon referred to as the five reference years.

On the one hand, applying \ref{mod:2_timeseries} to the \emph{gurobipy} model leads to robust solutions with the lowest cost increase on average. As can be seen in Figure \ref{fig:mod2_Results_gurobi}, the TAC is similar across all $40$ years. On average, additional cost of $5.6$bn€ or $10\%$ with a range of $[60.5, 63.15]$ are the result of this modification. The overall deviation from the average robust solution is in the range $[-1\%, 3\%]$ showing the effectiveness of the algorithm independent of the selected reference year.
The consistently lower costs indicate that \ref{mod:2_timeseries} is the best performing modification. 

\begin{figure}[h!]
\centering
  \includegraphics[width=\textwidth]{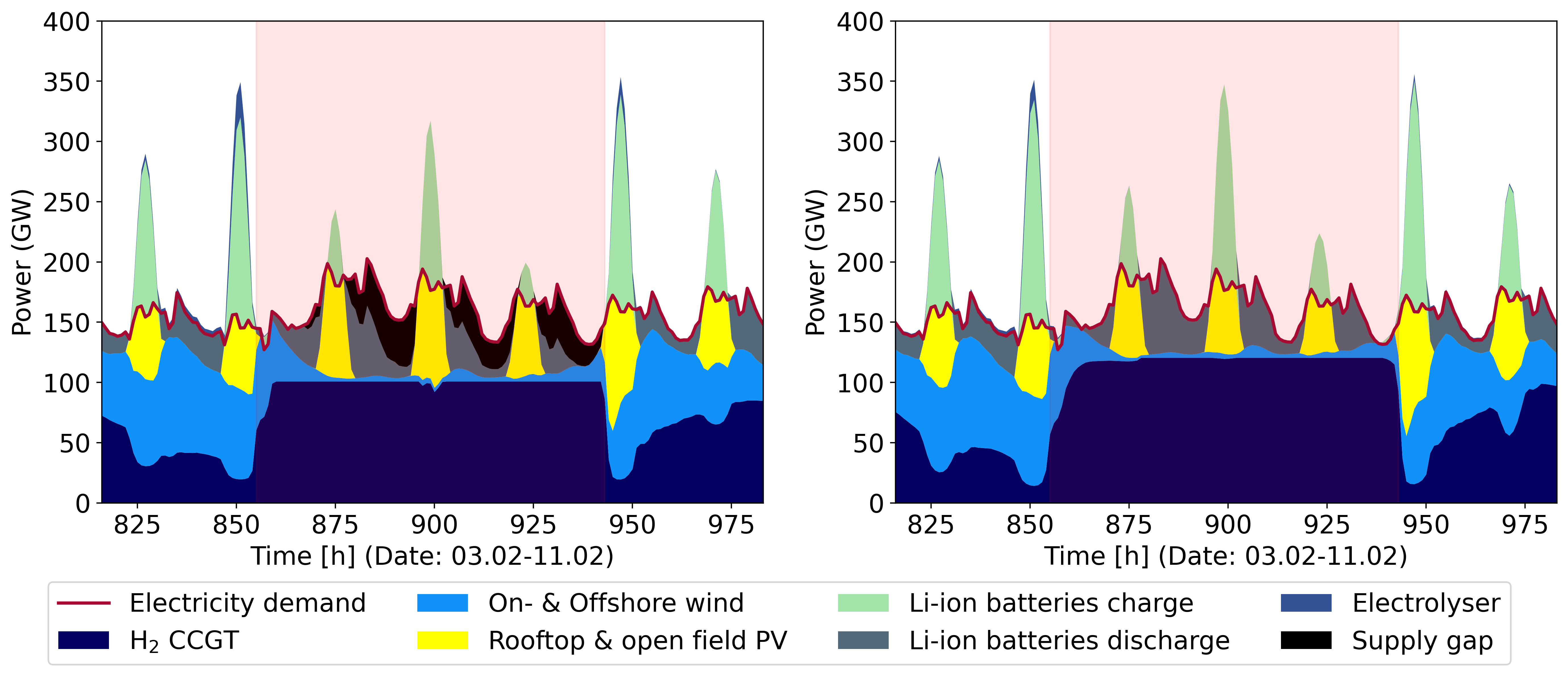}
  \caption{Feasibility testing of the energy system optimised for $1987$ in $1994$ before modification on the left and after modification on the right. The cold dark lull period in early February (see Figure \ref{fig:ExampleTimePeriods} f)) is marked in red. In the left graphic, due to insufficient backup capacity the supply gap variable has to be utilised meaning the energy system is not robust, i.e. there are still supply gaps after optimisation. In the right graphic, the energy is fully operational during the cold dark lull after applying \ref{mod:2_timeseries} to the original optimisation problem.}
  \label{fig:darklull}
\end{figure}

On the other hand, applying \ref{mod:2_timeseries} in ETHOS.FINE is shown in Figure \ref{fig:mod2_Results_FINE}. After modifying the five selected reference years with \ref{mod:2_timeseries} to make the solutions robust, the share of TAC for wind onshore decreases ($-2\%- 27\%$). 
Similarly, a decrease in total investment in transmission (electricity grid and hydrogen pipeline) is observed ($0\%- 25\%$).
A general increase in investment is seen for PV ($+4\%- +15\%$), Li--ion batteries ($-1\%- +26\%$) as well as the hydrogen sector ($+9\%- +18\%$) for the robust energy system designs. 
The increase of PV can be explained by its below average, but still relevant, availability during dark lulls combined with Li--ion batteries to cover daily fluctuations.
As visible in Figure \ref{fig:darklull}, PV is mainly utilised together with Li--ion batteries to cover the fluctuating part of the electricity demand, the CCGT cover the bulk of the electricity demand, while the generation from wind is negligible.
Hydrogen is utilised for electricity generation to a higher degree, since it can provide flexible additional energy supply, especially during dark lulls. The overall increase in cost compared to the average cost of each of the five reference years is $+12\%-13\%$, compared to the weather year $2018$ it is $10\%-12\%$ and compared to the most expensive single year, which is a lower (dual) bound on the objective, it is $2.9\%-5\%$. 
Thus, \ref{mod:2_timeseries} leads to robust and on average cheaper solutions in ETHOS.FINE similar to the results in the \emph{gurobipy} model.

This also implies that optimisation based on average or recommended reference years systematically underestimates the required cost for robust energy supply by $>10\%$ in ETHOS.FINE.

\begin{table*}[t]
\footnotesize
\begin{tabular}{l|ll|ll|ll|ll|ll}
    \hline
        Capacities & 1987 & 1987* & 1996 & 1996* & 2009 & 2009* & 2010 & 2010* & 2013 & 2013*\\ \hline
        Wind (onshore) $[GW]$ & 258 & 243 & 217 & 189 & 259 & 241 & 244 & 241 & 239 & 202\\ 
        Wind (offshore) $[GW]$ & 0.0 & 0.0 & 0.2 & 1.8 & 0.0 & 0.0 & 0.4 & 1.4 & 12.4 & 17.7\\
        PV (rooftop) $[GW]$ & 400 & 427 & 430 & 505 & 356 & 411 & 394 & 420 & 374 & 423\\ 
        PV (open field) $[GW]$ & 348 & 348 & 348 & 348 & 348 & 348 & 348 & 348 & 348 & 348\\
        Li--ion batteries $[GWh]$ & 722 & 840 & 742 & 888 & 684 & 770 & 770 & 765 & 635 & 802\\
        H$_2$ salt caverns $[TWh]$ & 195 & 211 & 191 & 234 & 202 & 237 & 204 & 217 & 169 & 206\\
        CCGT hydrogen gas $[GW]$ & 101 & 121 & 112 & 118 & 124 & 125 & 99 & 125 & 108 & 121\\
        Electrolysers $[GW]$ & 142 & 140 & 144 & 165 & 142 & 161 & 137 & 146 & 137 & 148\\ 
        Electricity grid $[GW]$ & 444 & 444 & 429 & 405 & 386 & 386 & 405 & 425 & 483 & 463\\ 
        Hydrogen pipelines $[GW]$ & 914 & 800 & 889 & 686 & 1029 & 686 & 800 & 686 & 686 & 571
\end{tabular}
\caption{Capacity results of optimising for $5$ selected years for Germany. Columns marked with a * indicate the  robust system designed with \ref{mod:2_timeseries}.}
\label{tab:Alg3_Results_capacities}
\end{table*}

Compared to the results of this work, Ryberg~\cite{Ryberg2020Diss} estimates a residual load of about $61$GW and additional backup capacity required of about $25$GW for Germany. The difference to the $118$GW--$125$GW found in this study can be attributed to the fact that in the integrated European setting that Ryberg~\cite{Ryberg2020Diss} used, dark lulls can be partly suppressed by electricity transmitted from regions not hit by that dark lull as well as differences in demand data. 

In either case, the cold dark lull period is the most critical for CCGT -- their installed capacity is mainly driven by a single dark lull period, as shown in Figure \ref{fig:darklull}.
Figure \ref{fig:darklull} also illustrates the effects of \ref{mod:2_timeseries}. The left graphic shows the result of the ESM optimised for 1987 when testing its feasibility in 1994 reveals a supply gap. After applying \ref{mod:2_timeseries}, the time period gets integrated into the optimisation problem and after reoptimising the supply can now be covered using existing capacities.

\paragraph{Modification \ref{mod:3_combine}} uses several types of cutting planes principles in one algorithm. It converges for all years, often only requiring one iteration of \ref{mod:3B_localbalance}. Sometimes, multiple iterations of \ref{mod:1_demand} are necessary as well. The effect of \ref{mod:3A_yearlybalance} is marginal: It does not effect model run times, nor results.

\begin{figure}[h!]
\centering
  \includegraphics[width=\textwidth]{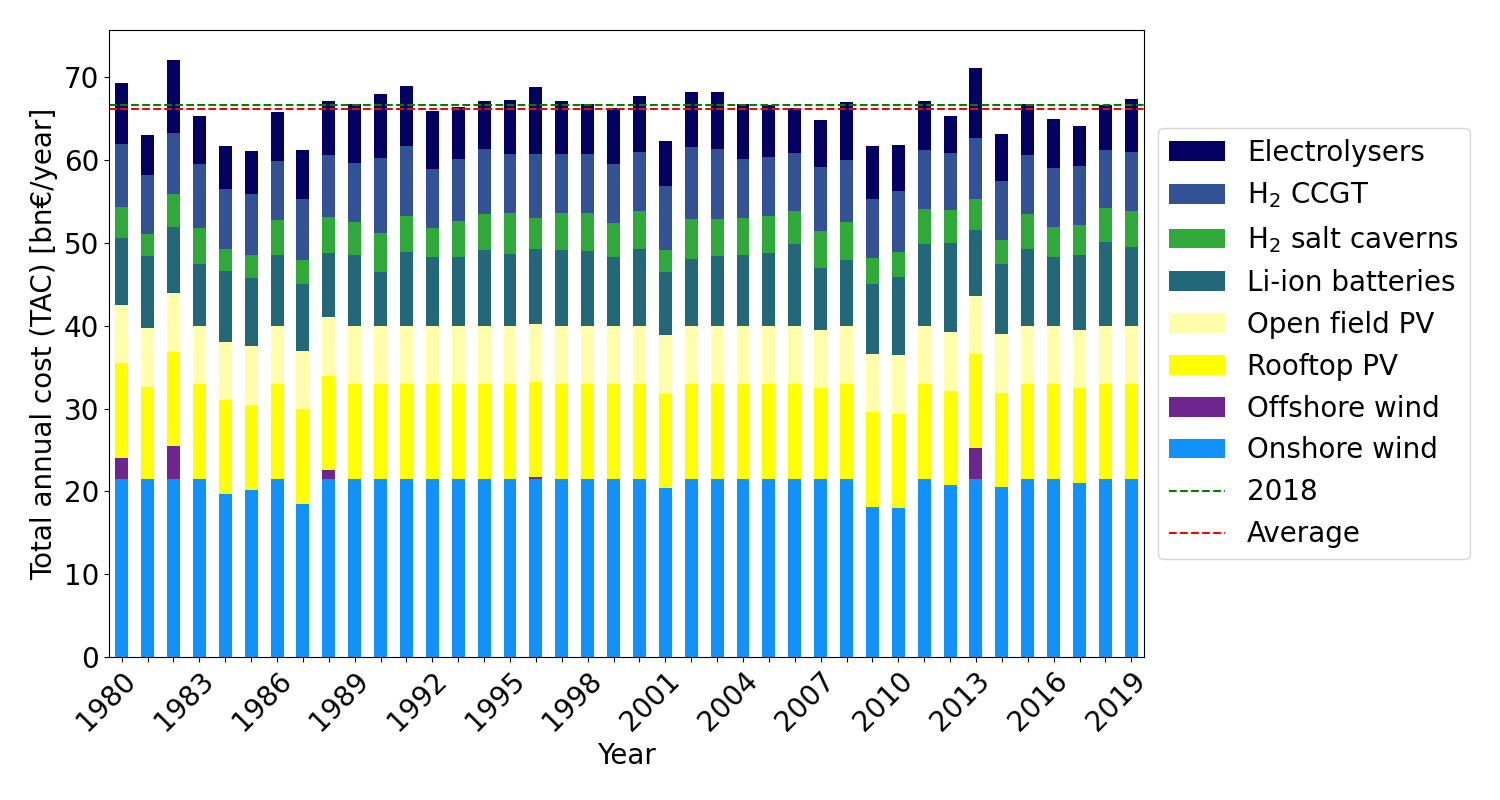}
  \caption{TAC comparison for robust solutions based on single years using \ref{mod:3_combine}, a single node \emph{gurobipy} model.}
  \label{fig:Alg8}
\end{figure}

The total costs average out to $66.0$bn€~per year, with a range of $[61.2,72.1]$, which is slightly more than \ref{mod:1_demand}. Thus, in terms of costs, the \ref{mod:1_demand} is preferable.
Figure \ref{fig:Alg8} gives the results of optimising each year independently for \ref{mod:3_combine}. Notably, the results for each year are very similar to each other, suggesting robust solutions share some traits. Cheap solutions contain less onshore wind.

\subsection*{Comparison of Modifications}
All three featured approaches generate robust solutions, i.e. solutions with no supply gaps, subject to numerical tolerances.
The modifications' performance is summarised in Table \ref{tab:all_algorithms}. 

\begin{table}[!ht]
    \footnotesize
    \begin{tabularx}{\textwidth}{|l|X|l|l|}
    \hline
        Modification & Increased Investment & Average & Least  \\ \hline
        \ref{mod:1_demand} & Li-Ion, $H_2$-CCGT, Salt caverns  & +5.7\% & +2.0\% \\
        Smoothed \ref{mod:1_demand} & $H_2$-CCGT, Salt caverns & +9.2\% & +2.4\% \\
        \ref{mod:2_timeseries} \emph{gurobipy} & PV, $H_2$-CCGT, Salt caverns, $H_2$-Electrolysers  & +2.5\%& +1.6\% \\
        \ref{mod:2_timeseries} ETHOS.FINE & PV, $H_2$-CCGT, Salt caverns, $H_2$-Electrolysers  & +3.7\%& +2.9\% \\
        \ref{mod:3_combine} & PV, Wind, $H_2$-CCGT, Salt caverns, $H_2$-Electrolysers & +10.6\% & +2.5\%  \\\hline
    \end{tabularx}
    \caption{Comparison of the different modifications in terms of convergence and performance. Average and least cost increase compared to best lower (dual) bound, the highest TAC of a single unmodified year.}
    \label{tab:all_algorithms}
\end{table}

The results show a moderate increase in investment for robust models leading to TAC increases of $1.6-2.9\%$ compared to the best available lower (dual) bound on the cost of a robust solution.
The dual bound is given by the maximum total cost for a single year, since the cost of a robust solution cannot be lower than the cost of a single year. 
Thus, the highest cost year is also the only one that could be robust\footnote{Clearly, this can be shown via a proof by contradiction, since each solution is assumed to be optimal and thus there must not be a feasible solution with a lower objective value.}.
For gurobipy, the lower bound is $59.7$bn€~annually. For the model in, ETHOS.FINE it is $113.6$bn€~annually.

A upper (primal) bound is given by the maximum capacity for each technology per region and year, since this will be feasible for all years. 
This results in a range of $[113.6,195.66]$bn€~annually for an optimal robust energy system in ETHOS.FINE.
For the \emph{gurobipy} model, the range for optimal robust solutions is $[59.7,67.6]$bn€~annually.

\paragraph{Capacity changes in near-optimal robust solutions}
It is important to acknowledge that any single robust solution may not accurately capture the properties of robust solutions as a whole, especially if considering the space of near-optimal solutions that may be relevant to decision makers, see the recommendations by Lombardi et al.~\cite{LOMBARDI2025}. 
While we did not explicitly perform \emph{modelling-to-generate-alternatives} on the solution space, computing $40$ extremal solutions on single years and iteratively robustifying them gives us a sample of $40$ near-optimal robust solutions that allow us to infer some properties of near-optimal robust solutions:

Compared to optimising individual years, near-optimal solutions to the robust model systematically use less onshore wind. 
This is plausible, since most years contain no extended (dark) lull periods coinciding with peak demands. 
In years without extended dark lulls, wind power provides stable and cheap energy, compared to PV that might require more storage and conversion units.
However, in years with dark lulls this advantage disappears.
As such, if costs for storage and conversion are priced in, optimising a year without a dark lull may lead to more investment in onshore wind than would be efficient.
Integrating appropriate dark lull periods, as suggested in this work, might help counteract that effect, leading to a more balanced energy mix.
A functioning capacity market, especially for backup technologies, is essential to ensure the needed capacities are installed and ready to generate electricity during dark lull periods. 
Overall, robust solutions were only $2-3\%$ more expensive compared to the lower (dual) bound given by the most expensive single year. 
Contrarily, a model based on average or recommended reference years systematically underestimates costs by over $10\%$.

\paragraph{Convergence}
The  \emph{gurobipy}  implementations of \ref{mod:1_demand} and \ref{mod:3_combine} tended to quickly converge. However, sometimes residual load differences in the order of MWh/kWh remained. In these cases, the algorithm was terminated after $20$ iterations.
\ref{mod:2_timeseries} needed between one and $54$ iterations to converge in the \emph{gurobipy} model. The number of iterations generally decreased for higher annual cost of individual years, but was mainly due to overwriting critical time periods where priority was given to critical time periods selected in earlier iteration steps. In ETHOS.FINE this issue was prevented by evaluating critical time periods on the modified data largely preventing overwriting. The calculations then needed between one and eight iterations for the five reference years. Here, one iteration incurs the same computational load of solving one ESM, or less if warm--starts help reducing computation times.

\subsection*{Full Load Hours and System Cost}\label{app:windplot}
\begin{figure}[h!]
    \centering
    \begin{minipage}[t]{0.49\textwidth}
    \includegraphics[width=\textwidth]{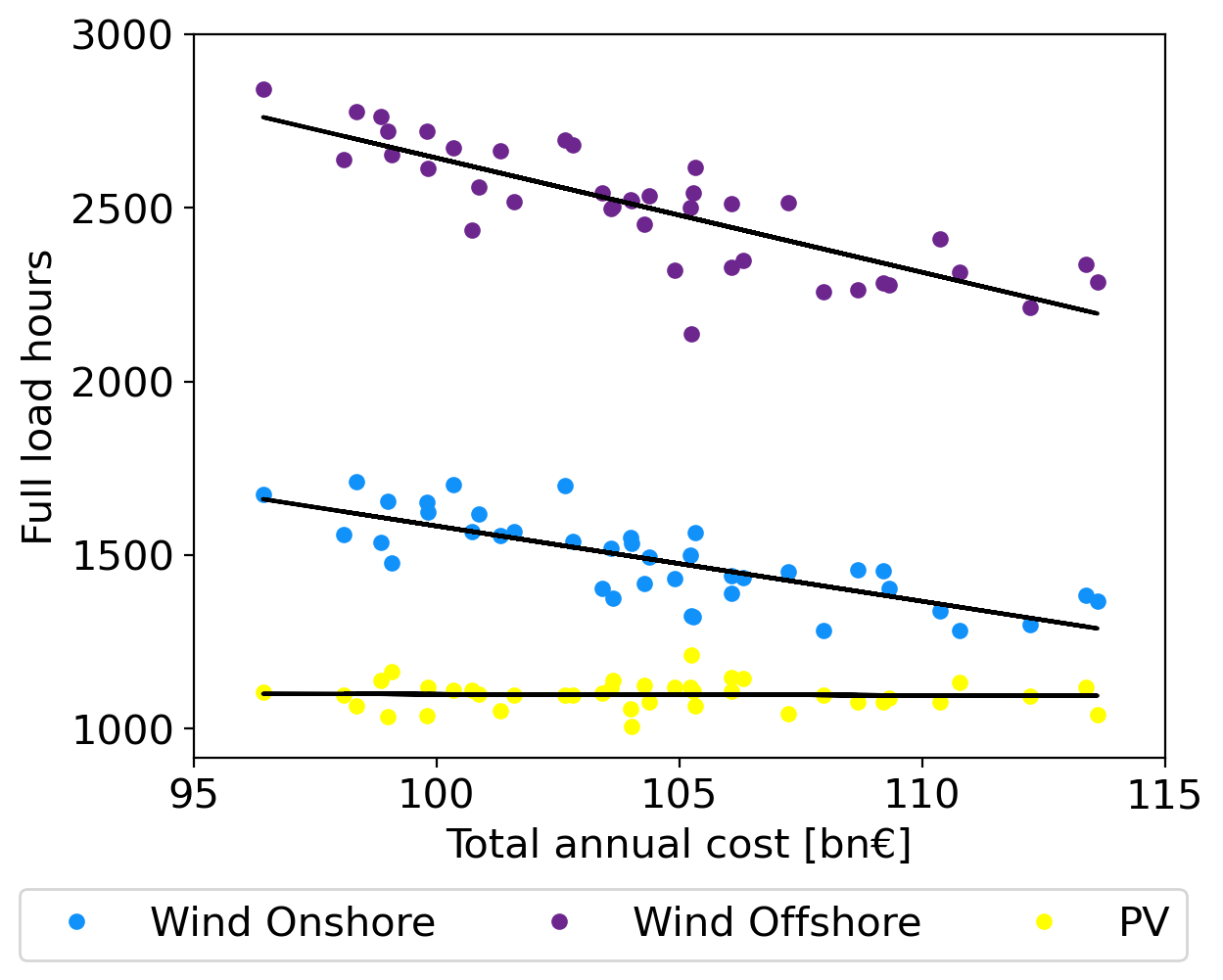}
    \end{minipage}
    \begin{minipage}[t]{0.49\textwidth}
        \includegraphics[width=\textwidth]{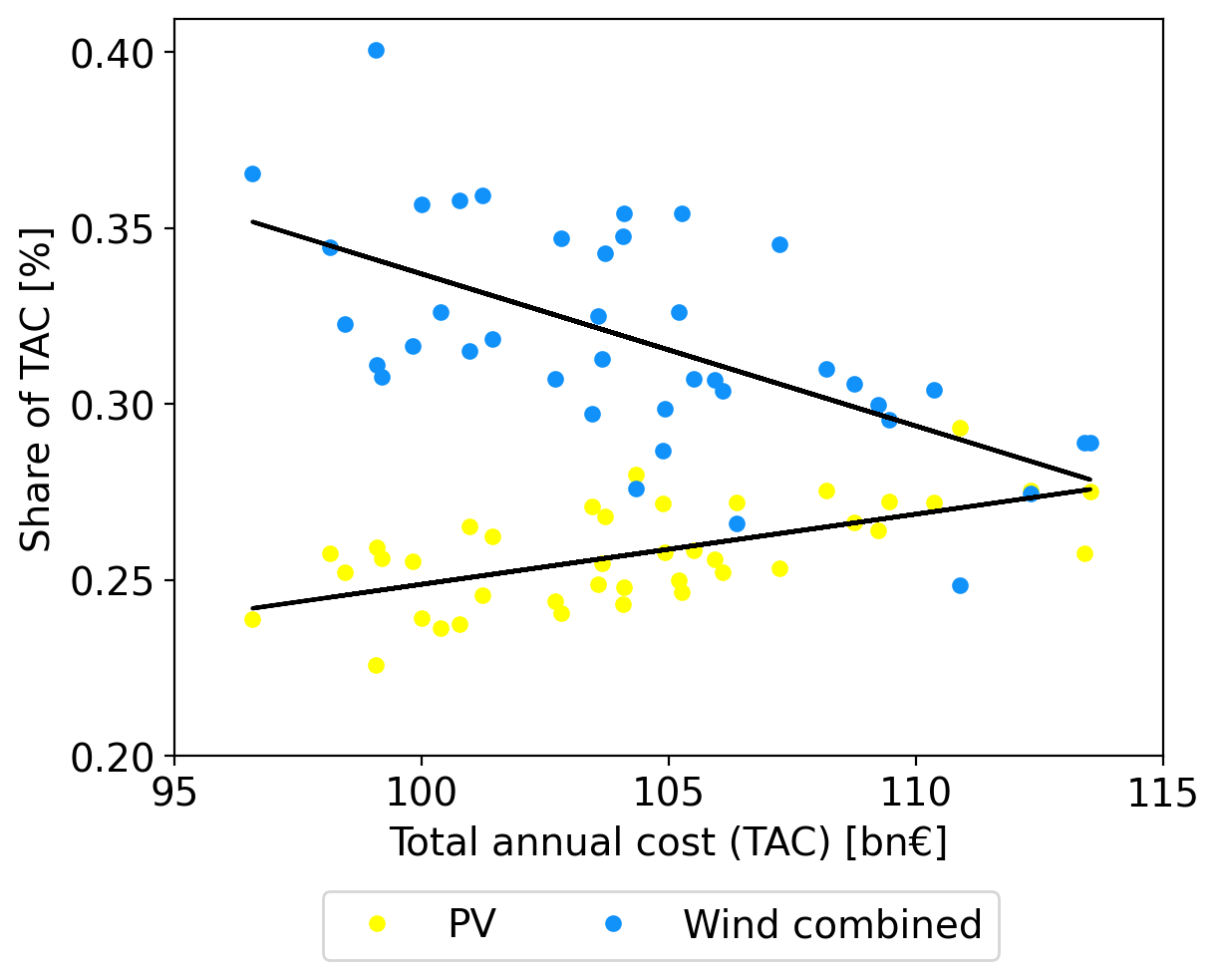}
    \end{minipage}
    \caption{The left graphic shows wind on-- and offshore as well as PV full load hours compared to TAC in for all years. 
    The right graphic shows combined cost for wind and for PV as share of the TAC compared to TAC for all years. Each dot represents one year. Lines indicate regression lines of the different energy generator types.}
    \label{fig:FLHvscostShareTACvsTAC}
\end{figure}

Lower investment in onshore wind capacity was a reoccurring pattern in the effects of modifications.
The left graphic in Figure \ref{fig:FLHvscostShareTACvsTAC} shows the annual full load hours (FLH) for wind on-- and offshore as well as PV compared to TAC for the respective models. For wind on-- and offshore these are strongly correlated (Pearson correlation coefficients of $-0.77$ and $-0.81$ for wind on-- and offshore, respectively). The FLH of PV and the TAC are nearly uncorrelated (Pearson correlation coefficient of $-0.04$). For wind this mirrors earlier results of Gotske et al.~\cite{gøtske2024designingsectorcoupledeuropeanenergy}, who showed similar correlations for an European System.

The right graphic in Figure \ref{fig:FLHvscostShareTACvsTAC} shows the share of TAC of PV and wind, combined for on-- and offshore, compared to to the TAC. For wind these are negatively correlated (Pearson correlation coefficient of $-0.59$) while for PV these are strongly positively correlated (Pearson correlation coefficient of $0.59$) indicating that as full load hours of wind drop, wind capacity gets replaced with PV capacity. 
The two diagrams in Figure \ref{fig:FLHvscostShareTACvsTAC} together show that the TACs depend strongly on the availability of wind.
In weather years with low full load hours for wind, an increase in PV capacities is observed indicating a higher reliance on PV in general which is also reflected in the robust years as investment increases for PV capacities for \ref{mod:2_timeseries} and \ref{mod:3_combine}. 

\subsection*{Application of Results}
The approaches outlined in this work can be used to make any ESM more robust to uncertain time--series data. As we showed, this serves to makes models more realistic and more suitable for real--world application preventing supply gaps for a moderate cost increase of $1.6$--$2.9\%$.

Since uncertainties in demands or costs, which appear in the right--hand side (demands) or the objective (cost) of a mixed-integer linear program (MILP) can be reformulated as constraint--wise uncertainty, the methods outlined in this work can be applied to cost uncertainty as well.
Another natural usage for the modifications proposed in this work is as part of a Benders decomposition framework. 
However, this was not the focus of this work, but the constraints that are added in each modification can equally well be used as feasibility/optimality cuts, as they serve to invalidate significant parts of the solution space. 
Here, the algorithmic performance for achieving robustness can be seen as a proxy for their potential value as Benders cuts.

\subsection*{Limitations of the study}

The main drawback in this case study is that only one time--series for electricity demand in $2050$ was used. The availability of future demand data is limited and is difficult to compare if taken from different sources. The selected data from Forschungsstelle für Energiewirtschaft e. V. (FfE)~\cite{ExtremosTertiarySector,ExtremosHouseholdSector,ExtremosTransportSector,ExtremosIndustrySector} was chosen since it includes a severe cold period rarely observed in Germany. In combination with the weather years the resulting operational conditions revealed several periods pairing low electricity generation with high demand as can be seen in Figure \ref{fig:ExampleTimePeriods}. Therefore a high degree of robustness can be assured. Further demand time--series might still provide additional insights. This is especially the case if using weather years with an overall cold winter time and therefore a high amount of heating degree days, and thus varying total and local demand.
The weather data used are based on historical weather years leaving out the effects of climate change, which very likely will influence the design of energy systems in the future and should be included in studies using this methodologies.
The model used to develop the presented method is a Germany model without the possibility of imports or exports. While this island model approach was necessary for method development, the interconnectedness of the European energy system, as well as (cold) dark lulls covering large parts of Europe, is likely crucial, but has been neglected in this study. The Germany model proved sufficient for model development, but a more general model is necessary for more holistic recommendations.
Assuming that years with disadvantageous distribution of sunny and wind hours as well as low full load hours are the exception, allowing the last time step to have a lower state of charge than the first would lower the conservatism of the system and reduce cost. 
This would require a measure of robustness that fully protects against a certain base uncertainty set, but allows using up some stored hydrogen for outlier events. One possible approach for that is outlined in Bärmann et al.~\cite{globalised}.
At the same time, the model only provides an operational schedule under the assumption of perfect foresight within one year. 
In a more complex model setting, computing an operational schedule dynamically throughout the year might lead to some efficiency losses in the usage of energy, requiring additional investment to counteract this. 
Perfect foresight is assumed for the optimization of energy systems. This provides the advantage of being able to plan capacities with full information about cold dark lulls and overall energy generation. In a real--world setting a higher degree of conservatism leading to additional generation or storage would be needed.

\subsection*{Summary of Implications for Planning Energy Systems}
Our work contributes to the mounting evidence that suggest that using a fixed reference year and planning a renewable energy system based on that is unsuitable for practice~\cite{schyska2021sensitivity,Ruggles2024ReliableWindSolar,PFENNINGER20171,haddeland2022effects,Ruhnau_2022,grochowicz2024using}.
We show that not only does this lead to misallocation of resources, but to a systematic underestimation of costs for and investment in storage and conversion, i.e. Li-Ion battery capacity and CCGT. Thus, if energy systems are planned based on single years, policy makers need to separately assess how much, not whether, additional storage and conversion capacities are needed and be aware of the fact that this will incur additional costs. 
For practitioners who model energy systems, we propose three workable approaches that can ensure energy system designs are robust against a range of weather realisations. We also note that the total amount of wind hours per year is strongly correlated with energy system costs. Thus, for a conservative estimate, low-wind reference years are better suited.

\section*{STAR Methods}\phantomsection\label{sec:methodology}

\subsection*{Key Resource Table}

\begin{table}[htbp]
\centering
\begin{tabularx}{\textwidth}{|X|X|}
\hline
\textbf{REAGENT or RESOURCE} & \textbf{SOURCE} \\
\hline
\multicolumn{2}{|l|}{\textbf{Deposited data}} \\
\hline
Load Curves of the Tertiary Sector – eXtremOS solidEU Scenario (Europe NUTS-3) &
\href{https://opendata.ffe.de/dataset/load-curves-of-the-tertiary-sector-extremos-solideu-scenario-europe-nuts-3/}{opendata.ffe.de / tertiary load curves} \\
\hline
Load Curves of the Household Sector – eXtremOS solidEU Scenario (Europe NUTS-3) &
\href{https://opendata.ffe.de/dataset/load-curves-of-the-private-household-sector-extremos-solideu-scenario-europe-nuts-3/}{opendata.ffe.de / household load curves} \\
\hline
Load Curves of the Transport Sector – eXtremOS solidEU Scenario (Europe NUTS-3) &
\href{https://opendata.ffe.de/dataset/load-curves-of-the-transport-sector-extremos-solideu-scenario-europe-nuts-3/}{opendata.ffe.de / transport load curves} \\
\hline
Load Curves of the Industry Sector – eXtremOS solidEU Scenario (Europe NUTS-3) &
\href{https://opendata.ffe.de/dataset/load-curves-of-the-industry-sector-extremos-solideu-scenario-europe-nuts-3/}{opendata.ffe.de / industry load curves} \\
\hline
40 years of wind and PV data in hourly resolution &
\href{https://www.renewables.ninja/}{renewables.ninja / wind PV data} \\
\hline
Tool for Renewable Energy Potentials – Database &
\href{https://zenodo.org/records/6414018}{zenodo.org / renewable potentials} \\
\hline
Council Regulation (EU) No 1059/2003 &
\href{https://eur-lex.europa.eu/legal-content/EN/TXT/?uri=CELEX:32016R2066}{eur-lex.europa.eu / EU regulation 1059} \\
\hline
\multicolumn{2}{|l|}{\textbf{Software and algorithms}} \\
\hline
Full code, model and data &
\href{https://github.com/FZJ-IEK3-VSA/Robust-Capacity-Expansion}{github.com / robust capacity expansion} \\
\hline
\end{tabularx}
\label{tab:key-resources}
\end{table}

\subsection*{Method Details}
After solving a initial CAPEX problem, the proposed algorithm consists of four main steps,  as shown in Figure \ref{fig:flowchart}: 1. feasibility testing using data of other years \footnote{In robust optimisation, \emph{scenarios} are commonly used as a term to denote different possible realisations of uncertainty. In this setting, each scenario consists of one year--long time--series of data for each type of renewable energy in combination with the demand affected by the weather. As such, this is a \emph{discrete uncertainty set}. While we introduce these notions for weather patterns only, the theory applies to any finite scenario set and their convex combination $\mathcal{U} = \textrm{conv}(u_1, ... u_n)$.}, 2. assessing loss of load/supply gaps, 3. identification of reasons for infeasibility, specifically critical time periods to derive modifications and, 4. solving the modified CAPEX problem before repeating steps 1-4.
Before describing these steps, we define robustness for ESMs.

\paragraph{Robustness in capacity expansion modelling}
\label{sec:problem_definition}
As noted before, general \ref{def:aro} has a bilevel structure: first stage decisions are made, some (worst--case) uncertainty realises, but then second stage decisions can be then made based on the then known uncertainty and the fixed first stage decisions.
For energy system modelling, we let the capacity decision variables be the first stage, since new power plants can not be build ad hoc to react to weather variations. 
The second stage decisions consist of operational decisions that can be made based on known weather variations. Thus, all second stage costs are operational expenditures (OPEX). We assume that the OPEX can be given as a flat percentage of investment costs/capital expenditure ($c_{CAPEX}^\top$), which is a reasonable assumption for many renewables such as wind power and solar PV~\cite{IRENA2017_capexopex}. This is assumed for all technologies included, but the fuel cost for hydrogen, which are endogenous, since hydrogen gas turbines can only use hydrogen produced by the system itself.

Via those simplifications, we remove costs from the second stage decisions. Mathematically, this is a crucial simplification as it allows us to collapse the ''min-max-min'' objective of general \ref{def:aro} into a simple minimisation problem.
A second motivation for doing that is that our approach emulates the principle of a Bender's decomposition (see the most recent book by Hooker~\cite{hooker2024} for an overview of Bender's decomposition). By removing the second stage costs, we remove the need for optimality cuts, which significantly simplifies computations. This results in the following, simplified \ref{def:aro} model:

\begin{align}
    \min_{x \in \mathcal{X}} \; & c_{CAPEX}^\top x& \label{def:aro-esm}\tag{ARO-ESM}\\
    \text{s.t.} \quad & A x + B(u) y(u) \geq b(u) &\forall u \in \mathcal{U}\notag \\
    & y(u) \in \mathcal{Y}(x,u) &\;\forall u \in \mathcal{U}.\notag
\end{align}

Here, we model supply and dispatch of electricity, (dis-)charge of storage systems as well as transmission and conversion of electricity and hydrogen for the technologies given in Table~\ref{tab:esm_components}. 

\begin{table}[h!]
\noindent
\centering
\begin{tabular}{l|l}
     & Technology \\ \hline
    Supply & 
    Rooftop PV, Open field PV,
    Onshore wind, Offshore wind\\
    Storage & 
    Li--ion batteries,
    H$_2$ salt caverns  \\
    Transmission & 
    Electricity grid,
    H$_2$ pipelines  \\
    Conversion & 
    Electrolysers,
    H$_2$ combined cycle gas turbines (CCGT)  \\
    Demand & 
    Electricity demand\\
\end{tabular}
\caption{Energy system components considered for development of the proposed methodology for optimising ESMs for Germany.}
\label{tab:esm_components}
\end{table}
    
Appendix \ref{app:lp_esm} gives an overview of the simplified model used for implementing constraints and variable demands to allow for an easy and fast testing of the multiple solution modifications used below. The full code, model and data are publicly available via GitHub~\cite{zenodo}.

\paragraph{Feasibility testing}
A solution to the mixed--integer program \ref{def:aro-esm} that is optimal and feasible for one year can lead to supply gaps for other years. 
Furthermore, since \ref{mod:2_timeseries} modifies the weather data resulting in artificial weather years, which are not based on atmospheric physics. The feasibility testing serves as a method to ensure that optimized energy systems can be operated in historic weather years.
Therefore, feasibility testing has two purposes: First, to detect supply gaps. Second, to ensure that optimized energy systems are realistic.
We call an ESM \emph{robust} against a set of years $\mathcal{U}$, if and if only it can be operated without supply gaps for any year $u \in \mathcal{U}$:

\begin{definition}[Robust energy systems]
    An energy system is robust against a set of time--series data $\mathcal U$ if and if only for each time--series, there exists a operation schedule that supplies all demand in time, while ensuring that the total amount of energy in storage is at least as big at the end of a year, as it was at the beginning.
\end{definition}

Remark that this definition explicitly defines robustness relative to a \emph{known} uncertainty set $\mathcal U$. Furthermore, it requires non--decreasing storage, which implies long--term supply security. That means even multiple sequential years with low full load hours for renewables will not deplete storage levels. This is a risk--averse strategy, as increased energy supply in beneficial years can not be used to offset lack of supply in bad years.

Now consider a solution $x_{u_i}, y(u_i)$ for one year $u_i \in \mathcal{U}$. To evaluate how well a given system would have performed in a different year ${u_j} \in \mathcal{U}, i \neq j$ it is necessary to determine whether for an optimized operating schedule any supply gaps remain. This is equivalent to reoptimising
\begin{align*}
    \min_{y(u_j)\in \mathcal{Y}(x_{u_i},u_j),\; \delta_{ij} \in \mathbb{R}^{\ge 0}} 1^\intercal\delta_{ij}\\ 
    s.t. \quad B(u_j)y(u_j) + \delta_{ij} \ge b-Ax_{u_i} \label{mip:comp}\tag{$COMP^i_{j}$} 
\end{align*}
Here, $\delta$ is a vector of supply gaps. It is implied that $\delta_{ij}$ is only added to constraints that cover energy supply/demand. 
Note that in the mixed integer program above, $x_{u_i}$, the capacity expansion, is a fixed input parameter retrieved from the selected year. 
Now an energy system designed for year $u_i$ is \emph{robust} for any year $u_j$ if and if only there exists an solution to the MILP $COMP_j^i$ above such that $\delta = 0$. If not, $\delta_{ij}$ encodes where supply gaps are and how large they are.\footnote{
Note that this is a relation between two $x_{s_i},x_{s_j}$ solutions, not their respective underlying time--series data, since there may be multiple (near-)optimal solutions for a given year. 
Computationally, this never proved to be an issue, so these notions are used interchangeably throughout this work. However, care should be taken if multiple capacity decisions have equal pricing.}

\paragraph{Optimisation problem modifications and clustering}
If an ESM displays supply gaps during feasibility testing, its CAPEX optimisation problem is modified to obtain a robust solution.
Consider a solution $\delta$ to \ref{mip:comp}. In this solution, $1^\intercal\delta_i^x$ units of energy are missing. This information needs to be integrated into the original CAPEX problem to enable the ESM to provide the extra energy supply required to meet demand.

Multiple modifications were iteratively developed to capture different sources of uncertainty.
Initial testing showed a subset of them to be ineffective. In order to also report null--results, these modifications are outlined in Appendix \ref{app:algorithms}.
The three promising modifications are inspired by three physical phenomena relevant to energy systems: 
1) The importance of peak demands for energy system design~\cite{Waite2020,Jackson2021}, 2) The existence of extended critical time periods~\cite{Ryberg2020Diss,Ruhnau_2022,grochowicz2023intersecting,grochowicz2024using} and 3) variations in total yearly energy supply and demand across different years~\cite{collins2018impacts}. 

For \ref{mod:1_demand}, update the demand vector $b_j$ of the starting year $j$ by adding the demand-supply gap $\delta$ to it, i.e. 
\begin{equation}
    b_j'= b_j+\delta. \tag{MOD 1 / Demand increase}\label{mod:1_demand}
\end{equation}
This adds the missing energy supply exactly when it is needed. One observation during initial testing was that this can lead to very large peaks in artificial demand, especially if done repeatedly over multiple iterations, leading to excessive battery installation. 
To counteract this, the artificial, additional demand can be divided up between neighbouring time periods, i.e. update $b$ and smooth $\delta$, i.e. $b'=b+f(\delta)$ for some function $f:\mathbb{R}^n_{\ge0} \mapsto \mathbb{R}^n_{\ge0} \setminus \emptyset$ that changes $b$ locally.
Figure \ref{fig:supply_gap} in Appendix \ref{app:algorithms} illustrates the effects of smoothing.
Throughout our work, we differentiate between smoothed and non--smoothed \ref{mod:1_demand}.

Since feasibility tests show that supply gaps usually occur in adjacent time steps, those are clustered to identify time periods with potential supply gaps, i.e. \textit{critical time periods}. Critical time periods are usually characterised by a combination of low electricity generation from renewables together with increased electricity demand, so--called \emph{(cold) dark lulls}. 
A hierarchical clustering was chosen, as this has been shown to be effective for time--series aggregation previously~\cite{tsam,bahl2016time}, and it allows for easy visual inspection and clustering along adjacent time--steps~\cite{hastie2009elements,clustering}.
For each time period, check whether the optimised capacities provide sufficient energy. 
For this, calculate the potential electricity generation of the selected ESM and compare this with the total demand in the time period.\footnote{Note that a positive average supply gap means that the investigated energy systems lack dispatchable capacity, and thus the feasibility testing of the respective time periods will always find a supply gap. However, a time period with negative average hourly supply gaps can still contain supply gaps that are only found during feasibility testing. This happens if either total annual electricity generation is too low to cover overall demand, or if there are some demand peaks that can not be covered.}
For calculating the potential electricity generation, full load is assumed for generation and backup capacities.
Let $b_{t,j}$ be the demand at hour $t \in T$, and $a_{t,j}^k\cdot x_s$ the time--series data and the calculated capacities of the energy system component $p \in P$ in year $u_j$, respectively. $P_{sup,con} = P_{sup} \cup P_{con}$ is used to denote that summation only happens over supply and conversion technologies.
Then, for time interval $T \subseteq \mathcal{T}$, the average hourly supply gap $\Delta(T)_j$ for some year $u_j'$ is given via

\begin{equation}
    \Delta(T)_j \coloneqq \quad \frac{\sum_{t \in T} \sum_{p \in P_{sup,con}} a_{t,j}^p\cdot x_s^p - b_{t,j}}{\sum_{t \in T}1}. \tag{$\Delta$(T)}\label{eq:delta_T}
\end{equation}

\ref{mod:2_timeseries}, iteratively inserts critical time periods into the original data to create synthetic time--series data for reoptimisation by first, selecting the ones with negative average hourly supply gap and second, when no critical time periods with negative hourly supply gap can be detected anymore, by reducing the sum of full load hours for PV and wind of the time period compared to the original data, i.e.
\begin{equation}
    a_{t,j}^p = a_{t,j'}^p \; \forall t \in T, p \in P \tag{MOD 2 / Synthetic time--series}\label{mod:2_timeseries}
\end{equation}
for some critical time period $T$ and years $j,j'$.

\ref{mod:3_combine} consists of multiple parts.
First, demand a weighted positive total energy balance for every year to ensure that there is sufficient energy overall.
Then, if during feasibility testing a positive slack $\delta^j$ is detected for year $j$, add
\begin{equation}
    \underbrace{\sum_{t\,\in\,\mathcal{T}}\left(\sum_{p\,\in\,P_{sup}} a^j_{tp}x_p\right)}_{\mathclap{\substack{\text{Renewable energy supply} \\ \text{in year $j$}}}}  -d(t) \ge \min\{0,1^\intercal\delta^j\}.
    \tag{MOD 3A / Yearly balance}\label{mod:3A_yearlybalance}
\end{equation}
to the model.
Here, the total demand is increased by the total absolute value of the slack vector, i.e.,  $1^\intercal\delta$. 
This is necessary since some energy will be lost during storage, conversion and transmission, which is not captured in the $a^j_{tp}$ model parameters and the supply capacities $P_{sup}$.

If this does not suffice to achieve feasibility for some year $j$, we model long-term energy storage more explicitly. In our model, $H_2$ is used for that.
The modelling is then done dynamically; for critical time periods $T$, auxiliary variables $\sigma_{j,T}$ are introduced that represent the use of CCGT power plants during periods of low renewable energy supply:
\begin{equation}
\underbrace{\sum_{t\,\in\,T}\left(\sum_{p\,\in\,P_{sup}} a^j_{tp}x_p\right)}_{\mathclap{\substack{\text{Renewable energy supply} \\ \text{in year $j$ during period $T$}}}} + \;\;\sigma_{j,T} -d(t) \ge 0.\tag{MOD 3B / Local $H_2$}\label{mod:3B_localbalance}
\end{equation}
At the beginning, the hydrogen storage is at $s^{H_2}_0$. In every hour after that, depending on the net energy balance, either more energy is stored or energy is taken from storage via $H_2$ to electricity conversion. At any point during a year, the hydrogen stored so far must be sufficient to at least cover all hydrogen demands during past/present critical periods. This can again be encoded as a MILP constraint:

\begin{equation}
\underbrace{s^{H_2}_0}_{\mathclap{\substack{\text{Initial $H_2$} \\ \text{storage}}}}+ \; \alpha \underbrace{\sum_{\substack{t\,\in\,\mathcal{T},\\\,t\,<\,T}}\left(\sum_{p\,\in\,P_{sup}} a^j_{tp}x_p\right) -d(t)}_{\mathclap{\substack{\text{Net renewable energy balance} \\ \text{in year $j$ to up time period $T$}}}} \;\ge \underbrace{\sum_{\substack{T'\,\subseteq \mathcal{T}\\ \,T'\,<\,T}}\sigma_{T'}}_{\mathclap{\substack{\text{Minimum $H_2$ required for} \\ \text{year $j$ to up time period $T$}}}}.\tag{MOD 3 / Combine}\label{mod:3_combine}
\end{equation}

The model parameter $\alpha$ encodes energy losses due to power to $H_2$ conversion. Note that this only is a lower bound on the energy required, since it is possible that more energy is used up during conversion in non--critical time periods or due to storage or transmission losses.\\

After reoptimisation only small supply gaps should remain. Subsequently \ref{mod:1_demand} that introduces additional demand is used to ensure leftover short-term supply gaps are covered, i.e the remaining supply gaps are added as an artificial electricity demand in the time periods where they appear. 

\paragraph{Case study of a German ESM}\label{sec:implementation}
The algorithm and all modifications were implemented and evaluated based on a fully renewable ESM for Germany in $2045$ that includes all technologies listed in Table \ref{tab:esm_components}.

Two implementations were evaluated. First, a benchmark model was implemented in \emph{gurobipy} that allows for comparison of all algorithms. A description of this model can be found in Appendix \ref{app:lp_esm}
Then, a more detailed $38$ node model was implemented in the open source Framework for Integrated Energy System Assessment (ETHOS.FINE) ~\cite{welder2018spatio} within the ETHOS modelling suite~\cite{gross2023ethos} to ensure validity of the results. This includes time--series aggregation via the TSAM package~\cite{tsam}.
Since ETHOS.FINE currently does allow interacting with the solver during the solution process, which is necessary to add cutting planes, this was only done for \ref{mod:2_timeseries}.

\paragraph{Case study data} $40$ years of wind and PV data in hourly resolution were taken from \href{renewables.ninja}{renewables.ninja}, using data originally provided by Staffell and Pfenninger~\cite{pfenninger2016long,staffell2016using}. They provide decades of hourly open data for both wind and PV allowing for easy integration into any kind of model.
The regional maximum capacity potentials for wind and PV were taken from Risch et al.~\cite{risch_2022_6414018, RischPotentials2022}, as this is highly detailed and validated data for Germany that if freely available. 
The basis for geodata is the \emph{Nomenclature des Unités territoriales statistiques} –- NUTS, a classification of the European Union. Level $2$ of this classification~\cite{NUTS2EU} is used in this work.
Electricity demand data in hourly resolution for the year $2050$ for tertiary, household, transport and industry sectors were taken from the \emph{Forschungsstelle für Energiewirtschaft e. V. (FfE)}~\cite{ExtremosTertiarySector,ExtremosHouseholdSector,ExtremosTransportSector,ExtremosIndustrySector}. 
On one hand, this data was selected because it covers a broad range of future electricity demands including cooling and heating demands. 
On the other hand, the data was calculated using the weather year $2012$ as a basis. Among the years considered in that study, 2012 was a year with an average amount of heating degree days, which fares well with our assumptions that storage levels at the first and last time steps need to be equal. The rare and extreme cold period with temperatures of $-30^{\circ}$C occurring in early February of $2012$ leads to a spike in heating demand making it ideal as a basis for highly robust ESMs. 
In combination with the dark lull in the weather data from $1994$, this represents a severe, although short ($< 4$ days) and extremely rare event in Germany and a high level of robustness for the results is to be expected. According to the German weather service~\cite{DWD1929}, the lowest temperature since the beginning of recordings is $-37.8^{\circ}$C. 
The utilised data and the developed code for either of the following calculations are publicly available via GitHub~\cite{zenodo}. The cluster specifications are outlined in Table \ref{tab:CAESAR}.

\begin{table}
    \centering
    \begin{tabular}{l|l}
        \hline
        CPU  & Intel(R) Xeon(R) Gold 6334\\
        Cores per node & 16\\
        Threads per node & 32\\
        Threads used & up to 3\\
        CPU max frequency & 3.6GHz\\
        RAM & up to 50GB (2TB available)
    \end{tabular}
    \caption{CAESAR computing cluster specifications of utilised nodes for optimising the ESM for Germany.}
    \label{tab:CAESAR}
\end{table}

\newpage

\section*{QUANTIFICATION AND STATISTICAL ANALYSIS}

There are no quantification or statistical analyses to include in this study.

\section*{RESOURCE AVAILABILITY}


\subsection*{Lead contact}


Requests for further information and resources should be directed to and will be fulfilled by the lead contact, Sebastian Kebrich (s.kebrich@fz-juelich.de).

\subsection*{Materials availability}
This study did not generate new materials.

\subsection*{Data and code availability}
The full code, model and data are publicly available via GitHub~\cite{zenodo}. Any additional information required to recreate the results reported in this paper is available from the lead contact upon request. 

\section*{ACKNOWLEDGMENTS}
This work is co-funded by the Deutsche Forschungsgemeinschaft (DFG, German Research Foundation) – 2236/2.
This work was partly supported by the Helmholtz Association as part of the Platform for the Design of a Robust Energy System and Raw Material Supply (RESUR) and the program, “Energy System Design.” 
This work was partly funded by the European Union (ERC, MATERIALIZE, 101076649). Views and opinions expressed are however those of the authors only and do not necessarily reflect those of European Union or the European Research Council Executive Agency. Neither the European Union nor the granting authority can be held responsible for them.

\section*{AUTHOR CONTRIBUTIONS}


\textbf{Sebastian Kebrich:} Conceptualisation, Methodology, Software, Validation, Formal analysis, Investigation, Data Curation, Writing -- Original draft, Writing -- Review \& Editing, Visualisation\\
\textbf{Felix Engelhardt:} Conceptualisation, Methodology, Software, Validation, Formal analysis, Investigation, Writing -- Original draft, Writing -- Review \& Editing, Visualisation\\
\textbf{David Franzmann:} Supervision, Writing -- Original draft, Writing -- Review \& Editing \\
\textbf{Christina Büsing:} Funding acquistion, Resources, Supervision, Writing -- Review \& Editing\\
\textbf{Heidi Heinrichs:} Funding acquistion, Resources, Supervision,
Writing -- Original draft, Writing -- Review \& Editing\\
\textbf{Jochen Linßen:} Resources

\section*{DECLARATION OF INTERESTS}
The authors declare no competing interests.

\bibliography{references}

\bigskip


\appendix

\section*{Appendix}
\section{Reference Linear Program}\label{app:lp_esm}

The following models gives a high--level explanation of the ESM we consider. For a detailed model, we refer to the implementation provided via GitHub~\cite{zenodo}. 
Consider a single node model. Let $t \in T$ denote time and $p \in P = P^{prod} \cup P^{stor} \cup P^{conv}$ denote different technologies, i.e. production, conversion and storage. Then, define the following variables:
\begin{align*}
    x_p \in [0,P^{max}_p] \subseteq \mathbb{R} \;&-&\textrm{amount of technology $p$ built,}\\
    s_t^{El} \in [0,S^{max}_{el}] \subseteq \mathbb{R} \;&-&\textrm{electrical energy in storage at time $t$,}\\
    s_t^{H2} \in [0,S^{max}_{H_2}] \subseteq \mathbb{R} \;&-&\textrm{hydrogen in storage at time $t$,}\\
    \delta_t^{+} \in [0,C^{max}_{H_2}] \subseteq \mathbb{R}\;&-&\textrm{energy for $H_2$ conversion at time $t$,}\\
    \delta_t^{-} \in [0,C^{max}_{El}] \subseteq \mathbb{R} \;&-&\textrm{energy from $H_2$ conversion at time $t$,}\\
    \Delta_t \in [0,d_{t}] \subseteq \mathbb{R} \;&-&\textrm{load shedding at time $t$,}\\
\end{align*}
where the upper bounds are given by technical system specifications, and by the respective demand per time for the load shedding. 
For simplicity, electrical energy storage is modelled without conversion losses, given that those are comparatively small. Contrary to that, $H_2$ conversion and storage is modelled explicitly.
Furthermore, unit commitment is not explicitly modelled, instead, the modeled guarantees the possibility of sufficient production, which might require turning off some generators in practice.

Then, minimise over $$\sum_{p \in P}c_{p}x_{p}+M\sum_{t \in T} \Delta_t,$$ 
where $M\in \mathbb{R}^+$ is chosen sufficiently large so that load shedding is never used if an alternative is possible. 
Note that this objective also implies that operational costs can be modelled as a flat percentage of investment cost. 
The minimisation above is subject to the following constraints.

First, energy must be conserved at each day:
$$\sum_{p \in P^{prod}}a_{tp}x_p + \delta_t^{-} - \delta_t^{+} + s_{t-1}^{El} - s_{t}^{El} +\Delta_t\ge d(t) \; \forall t\in T.$$
Here, $a_{tp}$ are coefficients modelling the production of energy from technology $p \in P^{prod}$ at time time $t \in T$ that in general are dependent on the weather scenario chosen.
For ease of notation, $-1 \sim |T|$ is used for indexing time.

Second, hydrogen storage and conversion is modelled explicitly:
$$ s_t^{H_2} = s_{t+1}^{H_2} + a_j^{+}\delta_t^{+} - a_j^{-}\delta_t^{-}.$$
Here, the $a_j^{+/-}$ coefficients are solely dependent on the conversion technologies used.

Third, both electricity and hydrogen storage must be bounded be the storage capacity built:
$$ s_t^{H_2} \le x_{H_2,storage}, $$
$$ s_t^{El} \le x_{El,storage}. $$

Fourth, yearly net energy production must be zero over the course of a year:
$$ s_0^{H_2} = s_{|T|}^{H_2}, $$
$$ s_0^{El} = s_{|T|}^{El}. $$

\section{Techno--Economic parameters}
\label{app:TechnoEconomic}
Table \ref{tab:TechnoEconomic parameters} shows the techno--economic parameters used in this work.

\begin{table*}[t]
\centering
\def\arraystretch{1.3}
\begin{tabular}[width=\textwidth]{l|llll}
    Technology & CAPEX$_{2050}$ & OPEX$_{fix, 2050}$& Life time [a]& Source\\ \hline
    PV (Rooftop) & 474 $\frac{\textit{€}}{\textit{kW}}$ & 10$\frac{\textit{€}}{\textit{kW a}}$ & 20 &~\cite{tsiropoulos2018cost, Kelm2019Vorbereitung}\\
    PV (Open field) & 320 $\frac{\textit{€}}{\textit{kW}}$ & 5.4$\frac{\textit{€}}{\textit{kW a}}$ & 20 &~\cite{tsiropoulos2018cost}\\
    Wind (Onshore) & 1000 $\frac{\textit{€}}{\textit{kW}}$ & 25$\frac{\textit{€}}{\textit{kW a}}$ & 20 &~\cite{tsiropoulos2018cost, Kreidelmeyer2020Transform}\\
    Wind (Offshore) & 2530 $\frac{\textit{€}}{\textit{kW}}$ & 63$\frac{\textit{€}}{\textit{kW a}}$ & 20 &~\cite{robinius2020wege}\\
    Li--ion batteries & 131 $\frac{\textit{€}}{\textit{kWh}}$ & 3.3$\frac{\textit{€}}{\textit{kWh a}}$ & 15 &~\cite{stolten2022neue}\\
    H$_2$ salt caverns & 0.7 $\frac{\textit{€}}{\textit{kWh}}$ & 0.01$\frac{\textit{€}}{\textit{kWh a}}$ & 40 &~\cite{caglayan2020technical}\\
    Electricity grid & 0.86 $\frac{\textit{€}}{\textit{kW km}}$ & 0.03$\frac{\textit{€}}{\textit{kW km a}}$ & 40 &~\cite{etri2014energy}\\
    H$_2$ pipelines & 0.185 $\frac{\textit{€}}{\textit{kW km}}$ & 0.01$\frac{\textit{€}}{\textit{kW km a}}$ & 40 &~\cite{Caglayan_2021robust}\\
    Electrolysers & 350 $\frac{\textit{€}}{\textit{kW}}$ & 11$\frac{\textit{€}}{\textit{kW a}}$ & 10 &~\cite{stolten2022neue}\\
    CCGT hydrogen gas & 760 $\frac{\textit{€}}{\textit{kW}}$ & 23$\frac{\textit{€}}{\textit{kW a}}$ & 20 &~\cite{stolten2022neue}\\
\end{tabular}
\caption{Techno--economic parameters considered in this work.}
\label{tab:TechnoEconomic parameters}
\end{table*}

Additionally, in ETHOS.FINE we also assumed a self discharge of $0.004\%$ per hour for Li--ion batteries.

\section{Description of Modifications}
\label{app:algorithms}

\textbf{/ Add constraint to enforce security of supply for time periods with low supply and high demand}

\begin{equation}
    \tag{MOD4 / Local Constraints}\label{mod:4_local_constraints}
\end{equation}

Begin with an energy system based on a single year time--series data. After performing feasibility testing, consider the case where at least one time period $T \subseteq \mathcal{T}$ is found for which in some weather year $j$, the given energy system can not supply sufficient energy.

For the energy system to supply sufficient energy, the energy supply during that period needs to be at least as large as the demand. Therefore, add a constraint that
$$ \sum_{t\,\in\,T}\underbrace{\left(\sum_{p\,\in\,P_{sup}} a^j_{tp}x_p\right)}_{\mathclap{\substack{\text{Renewable energy supply} \\ \text{in year $j$ during hour $t$}}}} + \;\;\underbrace{a_{H_2}x_{H_2\,Gas\,CCGT}}_{\mathclap{\substack{\text{Potential energy supply via} \\ \text{CCGT during hour $t$}}}} -d(t) \ge 0.$$
Here, $x$ are the capacity variables and $d(t)$ is the energy demand at time $t$. Note that this contrary to the coefficients $a^j_{tp}$ that model changing weather, the CCGT' $a_{H_2}$ conversion parameter is weather-independent. 

Battery storage is also weather-independent, but generally has less capacity than gas caverns. Therefore, it is not included in the constraint. 
However, for a model with small time periods or large expected battery storage, those can be integrated analogously to the $H_2$ case by adding the maximum battery capacity to the left hand side of the constraint above.\\

\begin{figure}[h!]
    \centering
    \includegraphics[width=\textwidth]{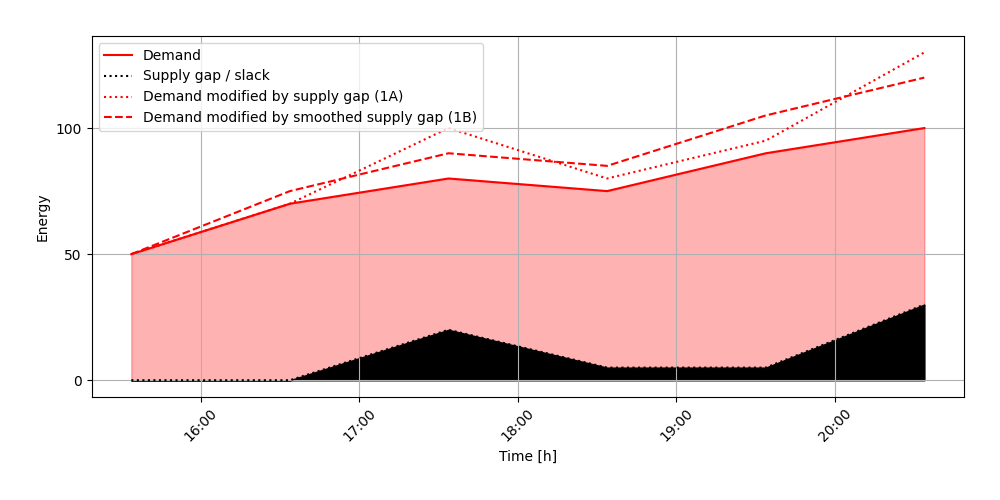} 
    \caption{Exemplary visualisation of the effect of smoothing for modification 1A and 1B on example data. The y-axis has no units, since this is for illustrative purposes only.} 
    \label{fig:supply_gap}
\end{figure}

\begin{equation}
    \tag{MOD6 / Global --$H_2$}\label{mod:6_global_H2}
\end{equation}

\textbf{/ Increase $H_2$ storage level at the end of year based on global supply gap} \\
Both of the approaches mentioned before are designed to compensated for supply gaps in specific time periods. In comparison, Modification $6$ aims to ensure sufficient energy supply for a full year as a whole.

For that, consider the vector of slack variables $\delta$ that encodes insufficient energy supply throughout the year.
We increase the required hydrogen storage level at the end of the year by the total absolute value of the slack vector, i.e.,  $1^\intercal\delta$. This ensures not only that more energy is supplied, but that the extra energy is supplied as $H_2$, which means it can be used  flexibly. \\

\paragraph{Comparison of modifications 1A and 1B}
Figure~\ref{fig:supply_gap} illustrates the differences between the two algorithms. For 1A, the supply gap is directly added to the demand time series. 
For 1B, smoothing of the supply gap lowers the peaks in the resulting demand profile, while introducing additional loads for nearby time periods.

\section{Evaluation of Ineffective Modifications}
None of the modifications that are outlined in Section~\ref{app:algorithms} are able to generate robust solutions on their own. Reasons for that are outlined below:

For Modification ~\ref{mod:4_local_constraints}, the algorithm will ensure that each cluster period is feasible, if sufficient hydrogen is available. 
However, the total available hydrogen is restricted by the total available energy and electrolyser capacity, which are not part of Modification \ref{mod:4_local_constraints}. Therefore, no convergence can be guaranteed in practice, due to insufficient total energy supply and $H_2$ conversion capacity.

Compensating this supply gap through the dynamic addition of hydrogen demands in Modification~\ref{mod:3B_localbalance} vastly improves performance. If gaps remain, they are in the range of hundreds of GWh, instead of tens of thousands of GWh in Modification ~\ref{mod:4_local_constraints}. However, this still does not lead to overall feasibility, if no mechanism for continuously adapting either the required extra supply or the total artificial energy demand is provided. 

For Modification~\ref{mod:6_global_H2}, the total available energy is increased eventually, but that does not imply that production capacities are sufficient during each time period.
Convergence is slow, often not achieving a robust solution even after $20$ iterations. Furthermore, costs tend to be high with large capacities for hydrogen production, but too few CCGT power plants to deal with extended dark lulls.

Finally, Modification~\ref{mod:3A_yearlybalance} suffers from non--convergence as well, even if large additional energy demands are added. To give an example, in one instance adding $>50\%$ to Germany's total energy demand, using the production parameters from another reference year, was insufficient to enforce feasibility for that specific weather year.

\end{document}